%% file: KinvDO.tex
\documentclass[a4paper,reqno,12pt]{amsart}
\usepackage[all]{xy}
\usepackage{bbm}
\usepackage{yhmath,mathrsfs}
\usepackage{eucal}

\input env.tex

\input ncpfaff_def.tex

\title%
[Principal symbols of $\Kc$-invariant differential operators]%
{On the principal symbols of $\Kc$-invariant 
 differential operators on Hermitian symmetric spaces}
\author{Takashi Hashimoto%
}
\address{
  Department of Informatics, 
  Graduate School of Engineering,
  Tottori University, 
  4-101, Koyama-Minami, Tottori, 680-8552, Japan
    }
\email{thashi@ike.tottori-u.ac.jp}
\date{\today}

\begin{document}

\begin{abstract}
Let $(G,K)$ be one of the following 
classical irreducible 
Hermitian symmetric pairs of noncompact type: 
$( \SU(p,q), S(\U(p) \times \U(q)) ),\;
( \Sp(n,\R), \U(n) )$,
or 
$( \SOstar(2n), \U(n) )$.
Let $\Gc$ and $\Kc$ be
complexifications of $G$ and $K$, respectively,
and let $P$ be a maximal parabolic subgroup of $\Gc$
whose Levi subgroup is $\Kc$.
Let $V$ be the holomorphic part of 
the complexifiaction of
the tangent space at the origin of $G/K$.
It is well known that
the ring of $\Kc$-invariant differential operators on $V$
has a generating system $\{ \varGamma_k \}$ 
given in terms of determinant or Pfaffian 
that plays an essential r{\^o}le in the Capelli identities
(\cite{HU91}).
Our main result of this paper is that
determinant or Pfaffian of the ``moment map'' 
on the holomorphic cotangent bundle of $\Gc/P$
provides a generating function for 
the principal symbols of $\varGamma_k$'s.
\end{abstract}

\maketitle

\section{Introduction}

Let $V:=\Alt_n$ be the vector space consisting 
of all alternating $n \times n$ complex matrices,
and $\C[V]$ the vector space consisting of 
all polynomial functions on $V$.
Then the complex general linear group $\GL{n}$ 
acts on $V$ by
\begin{equation}
\label{e:GL_action}
  g.Z:=g Z \tp{g}  \quad (g \in \GL{n}, Z \in V),
\end{equation}
from which 
one can define a representation $\pi$ of $\GL{n}$ 
on $\C[V]$ by
\begin{equation}
\label{e:GL_rep_on_polynom}
  \pi(g)f(Z):=f(g^{-1}.Z)   
  \quad  (g \in \GL{n}, f \in \C[V]).
\end{equation}
For $Z=(z_{i,j})_{i,j=1,\dots,n} \in V$,
with $z_{j,i}=-z_{i,j}$,
let $M:=(z_{i,j})_{i,j}$ and $D:=({\pd}_{i,j})_{i,j}$
be the alternating $n \times n$ matrices 
whose $(i,j)$-th entries are given by
the multiplication operator $z_{i,j}$
and the derivation $\pd_{i,j}:=\pd/\pd z_{i,j}$,
respectively.
Then the representation $\d\pi$ of $\gl_n$,  
the Lie algebra of $\GL{n}$,
induced from $\pi$ is given by
\begin{equation}
\label{e:gl_rep_on_polynom}
\d\pi(E_{i,j}) = -\sum_{k=1}^{n} z_{k,j} \pd_{k,i}
  \quad (i,j=1,2,\dots,n)
\end{equation}
where $E_{i,j}$ denotes the matrix unit of size $n \times n$
which is a basis for $\gl_n$.

Let us denote 
by $U(\gl_n)$ and $U(\gl_n)^{\GL{n}}$
the universal enveloping algebra of $\gl_n$,
and 
its subring consisting of $\GL{n}$-invariant elements,
respectively.
Let us denote 
by $\mathscr{PD}(V)$ and $\mathscr{PD}(V)^{\GL{n}}$
the ring of differential operators on $V$ 
with polynomial coefficients,
and its subring consisting of 
$\GL{n}$-invariant differential operators,
respectively. 
Then the following fact is known:
\begin{thm*}
[\cite{HU91}]
(1)\; 
 The ring homomorphism 
 $U(\gl_n)^{\GL{n}} \rightarrow \mathscr{PD}(V)^{\GL{n}}$
 induced canonically from $d\pi$ is surjective.
 
(2)\;
 For $k=0,1,\dots,\floor{{n}/{2}}%
      \footnote{%
        For $x \in \R$, 
        $\floor{x}$ stands for the greatest integer not exceeding $x$.%
      }$,
 let
 \begin{align}
  \varGamma_k:=\sum_{ {I \subset[n],\, |I|=2k} }
                        \Pf{z_I} \Pf{\pd_I},
     \label{e:skew_capelli}
 \end{align}
 where 
 the summation is taken over all subsets 
 $I \subset [n]:=\{ 1,2,\dots,n \}$ such that
 its cardinality is $2k$,
 and $z_{I}$, $\pd_I$ denote submatrices of $M,D$
 consisting of $z_{i,j}, \pd_{i,j}$ with $i,j \in I$.

 Then
 $\{ \varGamma_k\}_{k=0,1,\dots,\floor{n/2} }$ 
 forms a generating system for $\mathscr{PD}(V)^{\GL{n}}$,
\end{thm*}
In particular, for $k=0,1,\dots,\floor{n/2}$,
there exist elements of $U(\gl_{n})^{\GL{n}}$
that correspond to $\varGamma_k$ under the homomorphism,
which are called \textit{skew Capelli elements}.
As the names show,
they play an essential r{\^o}le 
in the skew Capelli identity.

Now, following \cite{Itoh01,KW02},
let us consider an alternating $2n \times 2n$ matrix $\bs{\Phi}$ 
with entry in $\mathscr{PD}(V)$ 
given as follows:
\begin{equation}
\label{e:Phi}
\bs{\Phi}
:=\left[
\begin{array}{cccc|cccc} 
   0     & z_{1,2} & \cdots & z_{1,n} &     &      &   &   u              
    \\
-z_{1,2} & 0       & \ddots &  \vdots  &    &      & u &           
    \\  
\vdots   &\ddots   &  0     &z_{n-1,n} &    &\adots&   &
    \\[3pt]
-z_{1,n} &\cdots   &-z_{n-1,n}& 0      & u  &      &   &
    \\[3pt]  
  \hline
       &      &      & -u   & 0      &\pd_{n-1,n}& \cdots & \pd_{1,n}
    \\
       &      &\adots&      &-\pd_{n-1,n}& 0  &\ddots  &\vdots  
    \\
       &  -u  &      &      &\vdots  &\ddots  & 0  & \pd_{1,2} 
    \\[2pt]
  -u   &      &      &      &-\pd_{1,n}&\cdots & -\pd_{1,2}& 0
\end{array}    
\right],
\end{equation}
where $u \in \C$.
Though our original motivation of this work is 
to understand the skew Capelli elements more deeply,
let us focus on the corresponding 
\textit{commutative} objects 
i.e.~the principal symbols of $\varGamma_k$,
which we denote by $\gamma_k$ in this paper,
before we enter the \textit{noncommutative} world. 
So, it is immediate from the minor summation formula 
of Pfaffian (see \cite{IW06}, or \eqref{e:msf_pf_square} below)
that the principal symbol $\sigma(\Pf{\bs{\Phi}})$
of $\Pf{\bs{\Phi}}$ provides a generating function 
for $\{ \gamma_k \}$:
\begin{equation}
\label{e:gf_Pf_of_Phi}
 \sigma(\Pf{\bs{\Phi}}) 
  = \sum_{k=0}^{\floor{n/2}} u^{n-2k} \gamma_k.
\end{equation}
As for the noncommutative counterpart,
we can show that $\Pf{\bs{\Phi}}$ can be expanded 
in the same way as \eqref{e:gf_Pf_of_Phi},
but with the coefficient $u^{n-2k}$ of $\gamma_k$ replaced by
a certain monic polynomial 
in $u$ of degree $n-2k$ 
(see \cite{GFinvDO}). 
Therefore, 
if $\bs{\Phi}$ came from $U(\gl_n) \otimes \Mat{n}(\C)$,
it would follow immediately from the property of Pfaffian 
with noncommutative entry 
that $\varGamma_k$'s belong to $U(\gl_n)^{\GL{n}}$.
However, it is obvious from \eqref{e:gl_rep_on_polynom}
that there exist no elements of $U(\gl_n)$ 
that corresponds to the multiplication operator $z_{i,j}$,
nor the derivation $\pd_{i,j}$.

What is the natural reason for considering 
the matrix $\bs{\Phi}$ above
(or, its commutative counterpart)?

We observe that the action \eqref{e:GL_action}
of $\GL{n}$ on $V=\Alt_n$ is 
the holomorphic part of the complexification of 
isotropy representation at the origin
of Hermitian symmetric space $\SOstar(2n)/\U(n)$.
Thus we embed $U(\gl_{n})^{\GL{n}}$ 
into $U(\so_{2n})^{\GL{n}}$
and seek for a generating function for $\{\gamma_k\}$
in the latter, 
in order to find an answer to the question raised above,
where $\so_{2n}$ denotes the complexification 
of the Lie algebra of $\SOstar(2n)$.

\begin{equation*}
\label{cd:capelli_map} 
\xymatrix{
  U(\gl_n)^{\GL{n}} \ar@{>>}[dd] \ar@{^{(}->} [dr]  &  \\
                        &  U(\so_{2n})^{\GL{n}} \ar[dl] \\
  \mathscr{PD}({V})^{\GL{n}}  &
 }
\end{equation*}

The real linear Lie group $\SOstar(2n)$
has irreducible unitary representations,
called the holomorphic discrete series representations.
Among them,
we consider a representation $\pi_{\lambda}$ 
constructed from a holomorphic character $\lambda$ 
of a maximal parabolic subgroup $P$
whose Levi subgroup is $\GL{n}$
via Borel-Weil theory.
Note that the representation space of $\pi_{\lambda}$
is a Hilbert space consisting of 
holomorphic and square-integrable functions 
defined on an open subset of $V$,
and that the restriction of $\pi_{\lambda}$ 
to $K$-finite part  
coincides with $\pi$ given in \eqref{e:GL_rep_on_polynom}
when $\lambda$ is trivial.

Let 
$\d \pi_{\lambda}$ be the differential representation 
induced from $\pi_{\lambda}$,
which we extend to the one of $\so_{2n}$ 
by linearity.
Take a basis $\{X_i\}$ for $\so_{2n}$,
and its dual basis $\{ X_i^{\vee} \}$,
i.e.~the basis for $\so_{2n}$ satisfying that
\[
 B(X_i,X_j^{\vee})=\delta_{i,j},
\]
where $B$ is the nondegenerate bilinear form 
on $\so_{2n}$ given by 
$B(X,Y):=\frac{1}{2}\trace{XY}$.
For $X \in \so_{2n}$ given, 
we denote by $\sigma(X)$ 
the principal symbol of
the differential operator $\d \pi_{\lambda}(X)$
substituted $\xi_{i,j}$ for $\pd_{i,j}$.

Now we define an element 
$
\sigma(\bs{X}) \in %
  \C[z_{i,j},\xi_{i,j};1 \leqsl i < j \leqsl n] \otimes \so_{2n}
$ 
by
\[
 \sigma(\bs{X}):=\sum_{i} \sigma(X_i^{\vee}) \otimes X_i.
\]
Note that $\sigma(\bs{X})$ is independent of 
the basis $\{ X_i \}$ chosen 
and that 
if $\lambda$ is trivial then
it can be considered the moment map 
on the holomorphic cotangent bundle $T^{*}(\SO_{2n}/P)$,
where $\SO_{2n}$ denotes the analytic subgroup of $\GL{2n}$
corresponding to $\so_{2n}$
(see \S 7).
We then rewrite $\sigma(\bs{X})$ 
using $\gamma_1$ and 
a newly introduced (commutative) indeterminate $u$,
which we denote by $\tild\sigma(\bs{X})$.
We show that
Pfaffian of $\tild\sigma(\bs{X})$ provides 
a generating function for $\{ \gamma_k \}$
(Corollary \ref{c:gen_fun_so}).
The reader will find that the matrix $\bs{\Phi}$
naturally appears on the way
(see Theorem \ref{t:moment_map_with_u}).

All the setup so far applies to 
the other two types of 
classical irreducible Hermitian symmetric pairs
of noncompact type
$(\SU(p,q), S(\U(p) \times \U(q)))$ 
and $(\Sp(n,\R), \U(n))$.
Let $(G,K)$ be one of the two.
Then 
the differential operators $\varGamma_k$ given above
have analogous objects, 
i.e.~$\Kc$-invariant differential operators 
acting on the space of polynomial functions 
on the holomorphic part of the complexification 
of the tangent space at the origin of $G/K$
which also play an essential r{\^o}le 
in the corresponding Capelli identity,
where $\Kc$ denotes a complexification of $K$.
In these cases,
they are given in terms of the sum of the product of
(minor-)determinants of matrices 
defined analogously to $M$ and $D$ given above
(see \eqref{e:Gammas} for their precise definitions).
If we define $\tild\sigma(\bs{X})$ 
similarly to the $\so_{2n}$ case,
we can show that
determinant of $\tild\sigma(\bs{X})$
provides a generating function for
their principal symbols
(Corollaries \ref{c:gen_fun_sp} 
and \ref{c:gen_fun_sl}).

The contents of this paper are as follows:
In Section 2,
we give realizations of $G=\SU(p,q), \Sp(n,\R)$ and $\SOstar(2n)$, 
and their complexification $\Gc$.
Then we make a brief review of construction of
the holomorphic discrete series representations
$\pi_{\lambda}$, 
via Borel-Weil theory,
which are induced from a character $\lambda$
of a maximal parabolic subgroup of $\Gc$
whose Levi subgroup is $\Kc$.
In Section 3,
we recall the definition of
the $\Kc$-invariant differential operators from \cite{IU01},
and introduce our main objects $\sigma(\bs{X})$ and
$\tild\sigma(\bs{X})$.
More concrete definitions of these objects
will be given in the subsequent sections case-by-case.
Then in Sections  4, 5 and 6,
according as $G=\SOstar(2n), \Sp(n,\R)$, 
or $\SU(p,q)$ $(p \geqsl q)$,
we explicitly calculate the differential operators
$\d \pi_{\lambda}(X_i)$
for a basis $\{ X_i\}$ of $\g:=\Lie(G) \otimes_{\R} \C$,
define $\sigma(\bs{X})$ and $\tild{\sigma}(\bs{X})$
explicitly,
and prove our main results that
Pfaffian or determinant of $\tild{\sigma}(\bs{X})$
provides a generating function for the principal symbol of
the $\Kc$-invariant differential operators. 
In the appendix,
we collect some minor summation formulae
of Pfaffian and determinant
which we make use of in proving the main results.

\section{Holomorphic Discrete Series}

Henceforth,
let $G$ denote
one of $\SU(p,q)\;(p \geqsl q)$, $\Sp(n,\R)$,
or $\SOstar(2n)$,
which we realize as follows:
\begin{align}
\label{e:realization_G}
 \SU(p,q) &=\{ g \in \SL{p+q}(\C); \tp{\bar{g}} I_{p.q} g=I_{p,q} \}, 
        \notag \\
 \Sp(n,\R) &=\{ g \in \SU(n,n); \tp{g} J_{n,n} g = J_{n,n} \},
               \\
 \SOstar(2n) &= \{ g \in \SU(n,n); \tp{g} J_{2n} g =J_{2n} \}.
        \notag
\end{align}
Here, for positive integers $p,q,n=1,2,\dots$, 
the matrices $I_{p,q}$, $J_{n}$, and $J_{n,n}$ are given by
\[
 I_{p,q}=\begin{bmatrix} 
          1_p &  \\  & -1_q 
         \end{bmatrix}, 
        \quad
 J_n=\begin{bmatrix}
      & & 1 \\[-2pt] & \adots & \\[-2pt] 1 & & 
     \end{bmatrix},
        \quad
 J_{n,n}=\begin{bmatrix}
         & J_n \\  -J_n & 
     \end{bmatrix},
\]
respectively.
Let $K$ be a maximal compact subgroup of $G$ given by
\(
 K:=\left\{ 
      \left[ \begin{smallmatrix} a & 0 \\ 0 & d \end{smallmatrix} \right] 
           \in G \right\} 
\),
where, for an element
\(
\left[ \begin{smallmatrix} a & 0 \\ 0 & d \end{smallmatrix} \right]
\)
in $K$, 
the submatrices $a$ and $d$ are of size
$p \times p$ and $q \times q$, respectively 
when $G=\SU(p,q)$,
or, are both of size $n \times n$
when $G=\SOstar(2n)$ and $\Sp(n,\R)$.
Let $\Gc$ and $\Kc$ be 
the complexifications of $G$ and $K$, respectively,
and $P$ the maximal parabolic subgroup of $\Gc$
given by
\(
 P:=\left\{ 
      \left[ \begin{smallmatrix} a & 0 \\ c & d \end{smallmatrix} \right] 
           \in \Gc \right\}
\)
so that its Levi subgroup equals $\Kc$.
Define a holomorphic character 
${\lambda}:P \to \C^{\times}$ by
\begin{equation}
\label{e:character_of_P}
 {\lambda}( \begin{bmatrix}
             a & 0 \\ c & d 
           \end{bmatrix}    )
 = (\det d)^{-s}
\end{equation}
for $s \in \Z$.
Denote by $\C_{\lambda}$
the one-dimensional $P$-module defined by 
$p.v=\lambda(p)v$ $(p \in P, v \in \C)$,
and by $L_{\lambda}$
the pull-back of the holomorphic line bundle
$\Gc \times_{P} \C_{\lambda}$
by the embbeding $G/K=G P/P \hookrightarrow \Gc/P$.

The space $\Gamma(L_{\lambda})$ 
of all holomorphic sections for $L_{\lambda}$
are identified with 
the space of all holomorphic functions $f$ on $G P$,
which is an open subset of $\Gc$,
that satisfy
\begin{equation}
\label{e:section} 
 f(xp)={\lambda}(p)^{-1} f(x) \quad (x \in G P, p \in P).
\end{equation}
Define $\Gamma^2(L_{\lambda})$ to be 
the subspace of $\Gamma(L_{\lambda})$ consisting of
square integrable holomorphic sections 
with respect to Haar measure on $G$:
\begin{equation}
\Gamma^2(L_{\lambda})
 :=\{ f \in \Gamma(L_{\lambda}); \int_{G} |f(g)|^2 \d g < \infty \},
\end{equation}
and an action $T_{\lambda}$ of $G$ 
on $\Gamma^2(L_{\lambda})$ by
\[
 (T_{\lambda}(g) f)(x):=f(g^{-1}x) \quad (g \in G, x \in G P).
\]
Then
$(T_{\lambda},\Gamma^2(L_{\lambda}))$ 
is an irreducible unitary representation of $G$,
called the \emph{holomorphic discrete series}
(if it is not zero).
One must impose some condition on $s$
in order to require that 
$\Gamma^2(L_{\lambda})$ be nonzero,
which we do not consider here
since we will be concerned with 
the representations of Lie algebras in this paper
(see e.g.~\cite{Knapp86} for the condition).

Now, according as
$G=\SU(p,q)\;(p \geqsl q),\Sp(n,\R)$, or $\SOstar(2n)$,
let us realize the bounded symmetric domain $\Omega$
as follows:

\begin{tabular}{ll}
 if $G=\SU(p,q)$, &
 $\Omega:=\{ z \in \Mat{p,q}(\C);1_q-\tp{\bar z}z>0 \}$;
    \\
 if $G=\Sp(n,\R)$, &
 $\Omega:=\{ z \in \Mat{n}(\C);1_n-\tp{\bar z}z>0, J_n \tp{z} J_n = z \}$;
    \\
 if $G=\SOstar(2n)$, &
 $\Omega:=\{ z \in \Mat{n}(\C);1_n-\tp{\bar z}z>0, J_n \tp{z} J_n = -z \}$,
\end{tabular}

\noindent
and let $G$ act on $\Omega$ by linear fractional transformation:
\begin{equation}
\label{e:linear_fractional_transformation}
 g.z=(az+b)(cz+d)^{-1} \quad 
    ( g=\left[\begin{matrix} 
               a & b \\ c & d 
              \end{matrix}
        \right] \in G, 
      z \in \Omega
    ).
\end{equation}
Then $\Omega$ is isomorphic to $G/K$ in each case.
If we denote by $\mathscr{O}(\Omega)$ 
the space of all holomorphic functions on $\Omega$,
define a map
\[
 \Phi: \Gamma(L_{\lambda}) \to \mathscr{O}(\Omega), \quad
    f \mapsto F
\]
by
\[
 F(z)=f( \left[
           \begin{matrix} 1 & z \\ 0 & 1 \end{matrix}
         \right] 
       ).
\]
Then $\Phi$ is a bijection.
Let $\H_{\lambda}:=\Phi(\Gamma^2(L_{\lambda}))$,
and define an action $\pi_{\lambda}$ of $G$ on $\H_{\lambda}$
so that the diagram \eqref{cd:T_and_pi} 
commutes for all $g \in G$:
\begin{equation}
\label{cd:T_and_pi}
 \xymatrix{ 
  \Gamma^2(L_{\lambda}) \ar[d]_{T_{\lambda}(g)} \ar[r]^{\Phi} 
                        &  \H_{\lambda} \ar[d]^{\pi_{\lambda}(g)} \\
  \Gamma^2(L_{\lambda}) \ar[r]^{\Phi} 
                        &  \H_{\lambda} \,. 
 }
\end{equation}
Explicitly, 
it is given by 
\begin{equation}
\label{e:pi_{lambda}}
 (\pi_{\lambda}(g)F) (z)= \det(cz+d)^{s} \, F( (az+b)(cz+d)^{-1} )
\end{equation}
for $g \in G$ and $F \in \H_{\lambda}$,
where 
$
g^{-1}=\left[ 
          \begin{smallmatrix} a & b \\ c & d \end{smallmatrix}
        \right]
$.

Let us introduce a few more notations.
Let $\g:=\Lie(G) \otimes_{\R} \C$, $\p:=\Lie(P)$,
$\bar \u$ the nilradical of $\p$,
and $\u$ its opposite. 
Then we have
$\g=\u \oplus \p$ and $\p=\mathfrak{k} \oplus \bar \u$,
where $\mathfrak k:=\Lie(K) \otimes_{\R} \C$.
Let $\d\pi_{\lambda}$ be
the differential representation of $\Lie(G)$
induced from $\pi_{\lambda}$,
which we extend to the one of 
the complex Lie algebra $\g$ by linearity.
Furthermore, 
identifying $\Omega$ with an open subset of $\u$
by 
$
z \leftrightarrow %
  \left[
    \begin{smallmatrix} 0 & z \\ 0 & 0 \end{smallmatrix}
  \right]
$,
let $U_{\Omega}$ be
the image of $\Omega \subset \u$ by the exponential map.

Since we realize the real linear Lie groups $G$ 
as in \eqref{e:realization_G},
the corresponding complex Lie algebras $\g$ are given by
\begin{align}
\label{e:realization_Lie(Gc)}
 \sl_{p+q} &= \{ X \in \Mat{p+q}(\C); \trace X=0 \}, 
     \notag   \\
 \sp_{n}  &= \{ X \in \Mat{2n}(\C); \tp{X} J_{n,n}+ X J_{n,n}=O \},
        \\
 \so_{2n}  &= \{ X \in \Mat{2n}(\C); \tp{X} J_{2n}+ X J_{2n}=O \},
     \notag
\end{align}
respectively.

\begin{remark}
For $X \in \g$,
if $t \in \R$ is sufficiently small 
then $\exp t X$ acts on $\Omega$ 
by \eqref{e:linear_fractional_transformation},
hence on $\H_{\lambda}$ by \eqref{e:pi_{lambda}},
so that its differential at $t=0$ coincides 
with $\d \pi_{\lambda}$.
\end{remark}

%
%

\section{Principal Symbols of $\Kc$-invariant Differential Operators}

Let $V$ denote the holomorphic part of 
the tangent space at the origin of $G/K$.
Then one can identify $V$ with $\u$ 
and construct a representation of $\Kc$ 
on the space $\C[V]$ of all polynomial function on $V$ 
through the action of $\Kc$ on $V=\u$.

Let $\mathscr{PD}(V)^{\Kc}$ 
be the ring of $\Kc$-invariant differential operators
with polynomial coefficient 
and let $r:=\Rrank G$, the real rank of $G$.
Then it is well known that
there exists a generating system 
$\{ \varGamma_k \}_{k=0,1,\dots,r}$ 
for $\mathscr{PD}(V)^{\Kc}$
which are given in terms of determinant or Pfaffian
as follows:
\begin{subequations}
\label{e:Gammas}
\begin{enumerate}
 \item 
 For $G=\SU(p,q)$ with $p \geqsl q$,
 \begin{equation}
 \label{e:Gamma_sl}
  \varGamma_k=\sum_{ \substack{I \subset[p], J \subset[q] \\ |I|=|J|=k} } 
                 \det(z^{I}_{J}) \det(\pd^{I}_{J}),
        \quad (k=0,1,\dots,q);
 \end{equation} 

 \item 
 For $G=\Sp(n,\R)$,
 \begin{equation}
 \label{e:Gamma_sp}
  \varGamma_k=\sum_{ \substack{I, J \subset[n] \\ |I|=|J|=k} }
                       \det(z^{I}_{J}) \det(\tilde{\pd}^{I}_{J}),
       \quad (k=0,1,\dots,n);
 \end{equation}

 \item 
 For $G=\SOstar(2n)$,
 \begin{equation}
 \label{e:Gamma_so}
  \varGamma_k=\sum_{ {I \subset[n], |I|=2k} }  
                \Pf{z_{I}} \Pf{\pd_{I}},
       \quad (k=0,1,\dots,\floor{n/2})
 \end{equation}

\end{enumerate}
\end{subequations}
(\cite{HU91}; see below for details).
We will find a generating function for 
the principal symbols of 
the differential operators $\varGamma_k$.

Define a $\Gc$-invariant nondegenerate bilinear form $B$ on $\g$ by
\begin{equation}
\label{e:bilinear_form}
B(X,Y)=
 \begin{cases}
 \; \trace{XY} \quad & \textrm{if }\g=\sl_{p+q},
              \\
 \; \frac12 \trace{XY} \quad & \textrm{if }\g=\so_{2n} \textrm{ or } \sp_{n}.
 \end{cases}
\end{equation}
Given a basis $\{ X_i \}_{i=1,\dots, \dim \g}$ for $\g$,
let us denote the dual basis with respect to $B$
by $\{ X_i^{\vee} \}$, 
i.e.~the basis for $\g$ satisfying
\[
 B(X_i,X_j^{\vee})=\delta_{i,j} 
\]
for $i,j=1,\dots,\dim \g$.
Then
we define an element $\bs{X}$ 
of $U(\g) \otimes \Mat{N}{(\C)}$ 
by
\begin{equation}
\label{e:ex_Casimir}
 \bs{X}:=\sum_{i=1}^{\dim \g} 
    X_{i}^{\vee} \otimes X_{i},
\end{equation}
where $U(\g)$ denotes the universal enveloping algebra of $\g$.
Namely,
we regard the former 
and the latter factors 
as elements of $U(\g)$ and $\Mat{N}{(\C)}$, 
respectively,
where
$N=p+q$ if $\g=\sl_{p+q}$,
and $N=2n$ if $\g=\sp_{n}$ or $\so_{2n}$.

Denoting by $\sigma(X)$ 
the principal symbol of $\d \pi_{\lambda}(X)$ 
for $X \in \g$,
let us define a $\g$-valued polynomial function
on the holomorphic cotangent bundle $T^{*}(\Gc/P)$
by 
\begin{equation}
\label{e:symbol_of_ex_Casimir} 
\sigma(\bs{X})
 :=\sum_{i=1}^{\dim \g}
    \sigma(X_{i}^{\vee}) \otimes X_{i}.
\end{equation}
By definition,
$\bs{X}$ and $\sigma(\bs{X})$ are independent of
the basis $\{ X_i \}$ chosen.

As we will explain below in more detail,
if the principal symbol $\sigma(X_i^{\vee})$ 
contains the parameter $s$,
we rewrite it using $\gamma_1$, 
the principal symbol of the Euler operator $\varGamma_1$ on $V$,
and introduce a new indeterminate $u$;
then we substitute $s-u$ into $\gamma_1 $in $\sigma(X_i^{\vee})$
which we denote by $\tild\sigma(X_i^{\vee})$.
If the principal symbol $\sigma(X_i^{\vee})$ does not contain $s$,
we set $\tild\sigma(X_i^{\vee}):=\sigma(X_i^{\vee})$. 
Finally, we define
\begin{equation}
\label{e:symbol_of_ex_Casimir_u} 
\tild\sigma(\bs{X})
 :=\sum_{i=1}^{\dim \g}
    \tild\sigma(X_{i}^{\vee}) \otimes X_{i}.
\end{equation} 
We will show that 
determinant or Pfaffian of $\tild\sigma(\bs{X})$
yields a generating function for the principal symbols 
$\{ \gamma_k \}_{k=0,1,\dots,r}$
of the generators $\{ \varGamma_k \}$ mentioned above 
in each case of (1), (2), and (3)
in the following sections.

\section{The Case $G=\SOstar(2n)$}

First let us consider the case 
where $G=\SOstar(2n)$, or $\g=\so_{2n}$.
In this section and the next,
$E_{i,j}$ denotes the $2n \times 2n$ matrix 
with its $(i,j)$-th entry being $1$ 
and all the others $0$.
Write $Z \in V=\u$ as
$Z=\sum_{i,j \in [n]} z_{i,j}E_{i,-j}$, $z_{j,i}=-z_{i,j}$,
where we agree that $-i$ stands for $2n+1-i$.
Let $M$ and $D$ denote alternating matrices of size $n$
whose $(i,j)$-th entry are given by
the multiplication operators $z_{i,j}$
and the differential operators $\pd_{i,j}:=\pd/\pd z_{i,j}$,
respectively.
Then
\[
  \varGamma_k:=\sum_{ {I \subset[n],\, |I|=2k} } 
                \Pf{z_{I}} \Pf{\pd_{I}}
       \quad (k=0,1,\dots,\floor{n/2})
\]
form a generating system for 
$\mathscr{PD}(V)^{\Kc}$ with $\Kc=\GL{n}$,
where
$z_{I},\pd_{I}$ denotes submatrices of $M,D$
consisting of the entries
whose row- and column- indices are 
both in $I \subset [n]$
(\cite{HU91}).

Take a basis $\{ X^{\epsilon}_{i,j} \}_{\epsilon=0.\pm;i,j \in [n]}$
for $\so_{2n}$ as follows:
{\allowdisplaybreaks
\begin{equation}
 \begin{aligned}
\label{e:basis_so}
 X^{0}_{i,j} &:= E_{i,j}-E_{-j,-i}  \quad ( 1 \leqsl i, j\leqsl n)
   \\
 X^{+}_{i,j} &:= E_{i,-j}-E_{j,-i} \quad (1 \leqsl i < j\leqsl n)
   \\
 X^{-}_{i,j} &:= E_{-j,i}-E_{-i,j} \quad (1 \leqsl i < j\leqsl n)
 \end{aligned}
\end{equation}
}

%
%
\begin{prop}
\label{p:rep_op_so}
The differential operators $\d \pi_{\lambda}(X^{\epsilon}_{i,j})$
$(\epsilon=0,\pm;i,j \in [n])$
are given by
{\allowdisplaybreaks
\begin{align}
 \d \pi_{\lambda}(X^{0}_{i,j})
      &= s \delta_{i,j} -\sum_{k \in [n]} z_{k,j} \partial_{k,i},
  \label{e:op_a_{i,j}}    \\
 \d \pi_{\lambda}(X^{+}_{i,j})
      &= -\partial_{i,j}, 
  \label{e:op_b_{i,j}}    \\
 \d \pi_{\lambda}(X^{-}_{i,j})
      &= -2s z_{i,j} 
         -\sum_{1 \leqsl k < l \leqsl n} 
             \left( z_{k,i}z_{j,l}-z_{k,j}z_{i,l} \right) \pd_{k,l}.
  \label{e:op_c_{i,j}}  
\end{align}
}
\end{prop}
\begin{proof}
Let us denote the $n \times n$ matrix 
whose $(i,j)$-th entry is $1$ 
and all the others are $0$ 
by $\umat{n}_{i,j}$,
and let $\tilde z=\sum_{k.l \in [n]} z_{k,l} \umat{n}_{k,l}$
with $z_{l,k}=-z_{k,l}$.
Then $z:=\tilde z J_n$ belongs to $\Omega$ 
if it is positive definite.

(I) First we calculate  $\d \pi_{\lambda}(X^{0}_{ij})$.
Writing 
$
\left[\begin{smallmatrix} a & 0 \\ 0 & d \end{smallmatrix}\right]
:= \exp (-t X^{0}_{ij})
$,
we see that
\begin{align*}
 \exp(-t X^{0}_{ij}).z 
  &= a z d^{-1}
     \\
  &= \left( 1_n -t \umat{n}_{i,j} + O(t^2) \right) z 
     \left( 1_n + t \umat{n}_{n+1-j,n+1-i}+O(t^2) \right)^{-1}
        \\
  &= \left( 
        \tilde z -t (\umat{n}_{i,j} \tilde z + \tilde z \umat{n}_{j,i}) 
     \right) J_n + O(t^2) 
        \\
  &= \Bigl(
       \sum_{k,l \in [n]}  z_{k,l} \umat{n}_{k,l}
       -t \sum_{k,l \in [n]} 
              z_{k,l} ( \umat{n}_{i,j} \umat{n}_{k,l} 
                       +\umat{n}_{k,l} \umat{n}_{j,i} ) 
     \Bigr) J_n + O(t^2)  
        \\
  &= \sum_{k,l \in [n]} 
       \left( z_{k,l}-t(\delta_{i,k} z_{j,l} + \delta_{i,l} z_{k,j}) \right) 
            \umat{n}_{k,l} J_n + O(t^2) 
\end{align*}
and that 
\[
 (\det d)^{s}
    =\left( \det(1_n + t \umat{n}_{n+1-j,n+1-i} + O(t^2)) \right)^s
    =e^{s t \trace{\umat{n}_{n+1-j,n+1-i}} + O(t^2)}.
\]
Therefore,
for $F \in \H_{\lambda}$,
we obtain that
\begin{align*}
  (\d \pi_{\lambda}(X^{0}_{ij})F)(z)
  &= \left. \frac{d}{dt} \right|_{t=0} 
           (\pi_{\lambda}(\exp(t X^{0}_{ij}))F)(z)
     \\
  &= \left. \frac{d}{dt} \right|_{t=0} 
      e^{st \trace{\umat{n}_{n+1-j,n+1-i}}}
       F \biggl( 
           \sum_{k,l \in [n]} 
             ( z_{k,l}-t(\delta_{i,k} z_{j,l} + \delta_{i,l} z_{k,j}) )
           \umat{n}_{k,l} J_n 
         \biggr)
     \\
  &= \biggl( 
       s \delta_{i,j} 
        - \sum_{k<l} (\delta_{i,k}z_{j,l}+\delta_{i,l}z_{k,j}) \pd_{k,l}
     \biggr) F(z) 
     \\
  &= \biggl( 
       s \delta_{i,j}
       - \sum_{i<l} z_{j,l} \pd_{i,l} 
       - \sum_{k<i} z_{k,j} \pd_{k,i} 
     \biggr) F(z)
     \\
  &= \biggl( s \delta_{i,j} -\sum_{k \in [n]} z_{j,k} \pd_{i,k} \biggr) F(z).
\end{align*}

(II) Next we calculate $\d \pi_{\lambda}(X^{+}_{ij})$.
Writing 
$
\left[\begin{smallmatrix} 1 & b \\ 0 & 1 \end{smallmatrix}\right]
 := \exp (-t X^{+}_{ij})
$,
we see that
\begin{align*}
 \exp(-t X^{+}_{ij}).z 
  &= z + b 
     \\
  &= ( \tilde z - t (\umat{n}_{i,j} - \umat{n}_{j,i}) ) J_n
     \\
  &= \sum_{k,l \in [n]}
      \left(
        z_{k,l}-t \delta_{i,k} \delta_{j,l} + t \delta_{j,k} \delta_{i,l}
      \right) \umat{n}_{k,l} J_n,
\end{align*}
and hence obtain that
\begin{align*}
   \d \pi_{\lambda}(X^{+}_{ij}) 
 &= - \sum_{k<l} 
      ( \delta_{i,k} \delta_{j,l} -\delta_{i,l} \delta_{j,k}) \pd_{k,l}
     \\
 &= - \pd_{i,j}.
\end{align*}  

(III) Finally, let us calculate $\d \pi_{\lambda}(X^{-}_{ij})$.
Writing 
$
\left[\begin{smallmatrix} 1 & 0 \\ c & 1 \end{smallmatrix}\right]
 := \exp (-t X^{-}_{ij})
$,
we see that
{\allowdisplaybreaks
\begin{align*}
 \exp(-t X^{-}_{ij}).z
  &= z(1+cz)^{-1}
     \\
  &=z \left( 1_n + t J_n (\umat{n}_{i,j}-\umat{n}_{j,i}) z \right)^{-1}
     \\
  &= \left( 
      \tilde{z} -t \tilde{z}(\umat{n}_{i,j}-\umat{n}_{j,i}) \tilde{z}
     \right) J_n + O(t^2)
     \\
  &= \Biggl(
     \sum_{k,l} z_{k,l} \umat{n}_{k,l}
    - t \sum_{a,b,c,d} 
         z_{a,b} z_{c,d} \umat{n}_{a,b}
          ( \umat{n}_{i,j}-\umat{n}_{j,i}) \umat{n}_{c,d} 
     \Biggr) J_n + O(t^2)
     \\
  &= \sum_{k,l} 
       \left( z_{k,l} 
             - t ( z_{k,i} z_{j,l}- z_{k,j} z_{i,l} )
       \right) \umat{n}_{k,l} J_n + O(t^2),
\end{align*}
}
and that
\begin{align*}
 \det(1+cz)^s 
  &=\left( 
      \det( 1_n + t J_n ( \umat{n}_{i,j}-\umat{n}_{j,i}) \tilde{z} J_n )
    \right)^s
  = e^{st \trace{( \umat{n}_{i,j}-\umat{n}_{j,i}) \tilde{z}} +O(t^2) }
     \\
  &= e^{ s t (z_{j,i}-z_{i,j}) + O(t^2) }
  =e^{-2 s t z_{i,j} + O(t^2) },
\end{align*}
from which we obtain that
\begin{align*}
 \d \pi_{\lambda}(X^{-}_{ij})
 &=-2s z_{i,j} 
   - \sum_{k<l} \left( z_{k,i}z_{j,l}-z_{k,j}z_{i,l} \right) \pd_{k,l}.
\end{align*}
This completes the proof.
\end{proof}
Noting that the dual basis of \eqref{e:basis_so}
is given by
\begin{align*}
 (X^{0}_{i,j})^{\vee}  &= X^{0}_{j,i},  
   &   
 (X^{\pm}_{i,j})^{\vee}  &= X^{\mp}_{i,j},
\end{align*}
let us define the element $\bs{X}$ 
of $U(\g)\otimes \Mat{2n}(\C)$ by \eqref{e:ex_Casimir};
it looks like
\begin{equation*}
\bs{X}
=\left[
\begin{array}{cccc|cccc}
X^{0}_{1,1} & X^{0}_{2,1} & \cdots & X^{0}_{n,1} 
         & X^{-}_{1,n} & \cdots & X^{-}_{1,2} & 0         
    \\
X^{0}_{1,2} & X^{0}_{2,2} & \cdots & X^{0}_{n,2} 
         & \vdots  & \adots & 0       & -X^{-}_{1,2}          
    \\  
    \vdots  & \vdots  &        & \vdots  
         & X^{-}_{n-1,n}& 0     & \adots   & \vdots
    \\[2pt]
X^{0}_{1,n} & X^{0}_{2,n} & \cdots & X^{0}_{n,n} 
         & 0       & -X^{-}_{n-1,n} &\cdots & -X^{-}_{1,n}
    \\[1pt]  \hline
X^{+}_{n,1} & \cdots  & X^{+}_{n-1,n}&  0    
         & -X^{0}_{n,n} & \cdots & -X^{0}_{n,2} & -X^{0}_{n,1}
    \\
\vdots  & \adots  & 0     & -X^{+}_{n-1,n} 
         & \vdots   &      & \vdots & \vdots 
    \\
X^{+}_{1,2} & 0       & \adots & \vdots  
         & -X^{0}_{2,n}  & \cdots & -X^{0}_{2,2} & -X^{0}_{2,1}
    \\[2pt]
  0     & -X^{+}_{1,2}& \cdots & -X^{+}_{1,n} 
         & -X^{0}_{1,n} & \cdots & -X^{0}_{1,2} & -X^{0}_{1,1}
 \end{array}
\right]. 
\end{equation*}

\begin{remark}
For the matrix $\bs{X}$ given above,
it is well known that Pfaffian%
     \footnote{Throughout the paper,
               for a given $2n \times 2n$ matrix $A$ 
               alternating along the antidiagonal,
               we denote $\Pf{A J_{2n}}$ by $\Pf{A}$ 
               for brevity.
     } 
$\Pf{\bs{X}}$ is a central element of 
the universal enveloping algebra $U(\so_{2n})$ 
(see \cite{IU01} or \cite{math.RT/0602055}
for definition and the properties 
of Pfaffian with noncommutative entry).
\end{remark}

Following the prescription 
\eqref{e:symbol_of_ex_Casimir},
let us define $\sigma(\bs{X})$ 
by substituting $\xi_{i,j}$ into $\partial_{i,j}$:
\begin{equation}
\label{e:symbol_of_X} 
\sigma(\bs{X})
 :=\sum_{\epsilon;i,j}
    \sigma( (X^{\epsilon}_{i,j})^{\vee}) \otimes X^{\epsilon}_{i,j}.
\end{equation}
%
%
%
%
\begin{thm}
\label{t:nilpotency_so}
Let
$u(z):=\exp \sum_{i<j} z_{i,j} X^{+}_{i,j} \in U_{\Omega}$.
Then we have
\begin{align}
 & \Ad(u(z)^{-1}) \sigma(\bs{X})
    = s \sum_{i} X^{0}_{i,i}
       - \sum_{i<j} \xi_{i,j} X^{-}_{i,j}
     \\
 &=\left[
   \begin{array}{cccc|cccc}
   s &               &           &     
       &\hph{-\xi_{1,n-1}} &\hph{-\xi_{2,n-1}} &\hph{\cdots}&\hph{-\xi_{2,n-1}}
    \\
     &        s     &            & 
        &               &           &            & 
    \\ 
     &              & \ddots     & 
        &               &           &            & 
    \\ 
     &              &            & s
        &               &           &             & 
    \\
   \hline
   -\xi_{1,n-1}  & -\xi_{2,n-1}   & \cdots & 0    
        & -s            &           &             &
    \\
   \vdots        & \vdots    & \adots        & \vdots 
        &               & \ddots    &             &
    \\
   -\xi_{1,2}    & 0          & \cdots     & \xi_{2,n-1}
        &           &           & -s          &
    \\
   0            &  \xi_{1,2} & \cdots     & \xi_{1,n-1}    
        &          &           &             & -s
   \end{array}
   \right].
\end{align} 
\end{thm}
\begin{proof}
This is just a simple matrix calculation,
but we give a rather detailed one,
which we will need in proving the theorem stated below.
In what follows,
for a matrix $A$ given,
let us denote the $(i,j)$-th entry of $A$
by $A_{ij}$.

Writing 
$
\left[ 
   \begin{smallmatrix} 1 & z \\ 0 & 1 \end{smallmatrix}
\right]
   :=u(z)
$
and
$
 \left[ 
   \begin{smallmatrix} A & B \\ C & D \end{smallmatrix}
 \right]
   := \sigma(\bs{X})
$,
we have 
\begin{equation}
\label{e:eq2show}
\Ad( u(z)^{-1} ) \sigma(\bs{X})
 = \begin{bmatrix}
    A-zC & Az-zCz+B-zD \\
    C    & Cz+D   
   \end{bmatrix}
\end{equation}

(I) First, we calculate the $(1,1)$ and $(2,2)$-blocks.
Let $\tilde{C}:=J_n C$.
Note that, by definition, 
$A_{ij}$ is $\sigma(X^{0}_{j,i})$. 
Since 
$(zC)_{ij}=(\tilde{z} \tilde{C})_{ij}$
equals
$\sum_{k=1}^{n} z_{i,k} \xi_{k,j}=-\sum_{k=1}^{n} z_{i,k} \xi_{j,k}$,
it follows from \eqref{e:op_a_{i,j}} that
\begin{equation}
\label{e:(1,1)-block_in_sigma}
  A = s 1_n + z C. 
\end{equation}
Then the fact that 
\eqref{e:eq2show} is an element of $\so_{2n}$ implies that
\begin{equation}
\label{e:(2,2)-block_in_sigma}
Cz+D=-J_n \tp{(A-zC)} J_n = -s 1_n.
\end{equation} 

(II) Next we show that the $(1,2)$-block equals $0$.
Let $\tilde{B}:=B J_n$.
Since 
$(\tilde{z} \tilde{C} \tilde{z})_{ij}$ equals
\[
 \sum_{k,l \in [n]} z_{i,k} \xi_{k,l} z_{l,j}
 =\Biggl( \sum_{k<l} + \sum_{k>l} \Biggr)\, z_{i,k} z_{l,j} \xi_{k,l}
 =\sum_{1 \leqsl k < l \leqsl n} 
    ( z_{k,i} z_{j,l} - z_{k,j} z_{i,l}) \xi_{k,l},
\]
it follows from \eqref{e:op_c_{i,j}} that 
$ \tilde{B} = -2s \tilde{z} -\tilde{z} \tilde{C} \tilde{z}$
and
\begin{equation}
\label{e:(1,2)-block_in_sigma}
 B = -2s {z} -{z} {C} {z}.
\end{equation}
Therefore we obtain that
\begin{align*}
 A z -z C z + B -z D 
  &=(s 1_n + z C)z - zCz - 2s z -zCz - z(-s 1_n -Cz)
     \\
  &=0. 
\end{align*}
This completes the proof.
\end{proof}

Now,
using the principal symbol 
$\gamma_1=\sum_{k<l} z_{k,l} \xi_{k,l}$ 
of the Euler operator $\varGamma_1$ on $V$, 
we rewrite 
$\sigma(X^{0}_{i,i})$ and $\sigma(X^{-}_{i,j})$,
which are the only principal symbols 
that contain the parameter $s$,
as follows:
\begin{align*}
 \sigma(X^{0}_{i,i})
  &= s - \gamma_1
          +\sum_{ \substack{ k<l \\ k,l \ne i} } 
                  z_{k,l} \xi_{k,l},
   \\
 \sigma(X^{-}_{i,j})
  &= -2s z_{i,j} + z_{i,j} \gamma_1 
     - \sum_{\substack{ k<l \\ \{k,l\} \cap \{i,j\} = \varnothing}}
        (z_{i,j}z_{k,l}-z_{k,i}z_{j,l}+z_{k,j}z_{i,l} ) \xi_{k,l}.
\end{align*} 
Then we substitute $s-u$ into $\gamma_1$ 
with a new indeterminate $u$ in the symbols above,
which we denote by
$\tild{\sigma}(X^{0}_{i,i})$ and $\tild{\sigma}(X^{-}_{i,j})$,
respectively:
\begin{align}
 \tild{\sigma}(X^{0}_{i,i}) 
  &:= u 
      +\sum_{ \substack{ k<l \\ k,l \ne i} } 
                  z_{k,l} \xi_{k,l},
   \label{tilde_sigma_X^{0}_{ii}}
     \\
 \tild{\sigma}(X^{-}_{i,j})
  &:= -(u+s) z_{i,j} 
     - \sum_{\substack{ k<l \\ \{k,l\} \cap \{i,j\} = \varnothing}}
        (z_{i,j}z_{k,l}-z_{i,k}z_{j,l}+z_{i,l}z_{j,k} ) \xi_{k,l}.
   \label{tilde_sigma_X^{-}_{ij}}
\end{align}
As for the other symbols,
set
$\tild{\sigma}(X^{\epsilon}_{i,j}):=\sigma(X^{\epsilon}_{i,j})$.
Then we define $\tild{\sigma}(\bs{X})$ by
\[
\tild{\sigma}(\bs{X})
 := \sum_{\epsilon;i,j} 
      \tild{\sigma}( (X^{\epsilon}_{i,j})^{\vee}) \otimes X^{\epsilon}_{i,j}.
\]

%
%
\begin{thm}
\label{t:moment_map_with_u}
Let $u(z) \in U_{\Omega}$ be as in Theorem \ref{t:nilpotency_so}.
Then we have
\begin{align}
 &\Ad(u(z)^{-1}) \tild{\sigma}(\bs{X})
  =(u+\gamma_1) \sum_{i} X^{0}_{i,i}
    - \sum_{i<j} \xi_{i,j} X^{-}_{i,j}
    - (s-(u+\gamma_1)) \sum_{i<j} z_{i,j} X^{+}_{i,j}
      \\
 &=\left[
   \begin{array}{cccc|cccc}
   u+\gamma_1 &               &           &     
       &-\tau\, z_{1,n-1}       & \cdots    &-\tau\, z_{1,2}    & 0
    \\
     & u+\gamma_1   &            & 
        &-\tau\,z_{2,n-1}      & \cdots    & 0          &\tau\,z_{1,2}
    \\ 
     &              & \ddots     & 
        & \vdots        & \adots    & \vdots     & \vdots
    \\ 
     &              &            & u+\gamma_1
        & 0             & \cdots    & \tau\,z_{2,n-2} & \tau\,z_{1,n-1}
    \\
   \hline
   -\xi_{1,n-1}  & -\xi_{2,n-1}   & \cdots & 0    
        & -u-\gamma_1   &           &             &
    \\
   \vdots        & \vdots    & \adots        & \vdots 
        &               & \ddots    &             &
    \\
   -\xi_{1,2}    & 0          & \cdots     & \xi_{2,n-1}
        &           &           & -u-\gamma_1 &
    \\
   0            &  \xi_{1,2} & \cdots     & \xi_{1,n-1}    
        &          &           &             & -u-\gamma_1
   \end{array}
   \right].
      \label{e:sigX2}
\end{align} 
Here, we set $\tau:=s-(u+\gamma_1)$ 
in \eqref{e:sigX2} for brevity.
\end{thm}
\begin{proof}
Let 
$
\left[ 
   \begin{smallmatrix} A(u) & B(u) \\ C & D(u) \end{smallmatrix}
\right]
  := \tild\sigma(\bs{X})
$.
(Note that the submatrix $C$ is the same as  
in the proof of Theorem \ref{t:nilpotency_so}).
Then it suffices to show that
\begin{gather}
 A(u) - z C = (u+\gamma_1) 1_n,
     \label{e:(1,1)-block_in_sigma_with_u}
     \\
 A(u)z - zCz + B(u) -zD(u) =(u+\gamma_1-s) z,
     \label{e:(1,2)-block_in_sigma_with_u}
     \\
 Cz + D(u) = -(u+\gamma_1) 1_n
     \label{e:(2,2)-block_in_sigma_with_u}
\end{gather}
(see \eqref{e:eq2show}).

(I) First we show \eqref{e:(1,1)-block_in_sigma_with_u} 
and \eqref{e:(2,2)-block_in_sigma_with_u}.
But the latter follows from the former,
as in the proof of Theorem \ref{t:nilpotency_so}.
Note that $A(u)_{ij}$ can be written as 
\begin{equation}
\label{e:A(u)_{ij}_so}
 A(u)_{ij}
 = (u + \sum_{k<l} z_{k,l} \xi_{k,l} ) \delta_{i,j}
    - \sum_{k=1}^{n} z_{k,i} \xi_{k,j}.
\end{equation}
Since the second summation in \eqref{e:A(u)_{ij}_so}
is equal to 
$-(\tilde{z}\tilde{C})_{ij}=-(zC)_{ij}$
as shown in (I) of the proof of Theorem \ref{t:nilpotency_so},
we obtain that 
\[
 A(u) = (u +\gamma_1) 1_n + zC.
\]

(II) Next we show \eqref{e:(1,2)-block_in_sigma_with_u}.
Let $\tilde{B}(u):=B(u) J_n$.
Since the summation part of 
$\tild\sigma(X^{-}_{i,j})=\tilde{B}(u)_{ij}$
equals 
$z_{i,j} \gamma_1 + (\tilde{z} \tilde{C} \tilde{z})_{i j}$
as in (II) of the proof of Theorem \ref{t:nilpotency_so},
we obtain that
\begin{equation}
\label{e:keyformula_(2,2)-block_in_sigma_with_u}
 B(u) = -(u+s+\gamma_1) z - zCz.
\end{equation} 
Now, combining \eqref{e:keyformula_(2,2)-block_in_sigma_with_u}
with \eqref{e:(1,1)-block_in_sigma_with_u} and
\eqref{e:(2,2)-block_in_sigma_with_u}, 
we obtain that 
\[
A(u)z-zCz+B(u)-zD(u)=(u+\gamma_1-s)z.
\]
This completes the proof.
\end{proof}

Note the similarity between 
the matrix $\bs{\Phi}$ given in \eqref{e:Phi} 
and $\tild{\sigma}(\bs{X})$ 
conjugated by $u(z)^{-1} \in U_{\Omega}$ 
given in \eqref{e:sigX2}.
It follows immediately from Theorem \ref{t:moment_map_with_u}
and the minor summation formula of Pfaffian \eqref{e:msf_pf_square}
that $\Pf{\tild{\sigma}(\bs{X})}$
yields a generating function for $\{ \gamma_k \}$.

%
%
\begin{cor}
\label{c:gen_fun_so}
Retain the notation above.
Then we have the following formula:
\begin{equation}
 \Pf{\tild{\sigma}(\bs{X})}
  =\sum_{k=0}^{\floor{n/2}} 
    (u+\gamma_1)^{n-2k} (s-(u+\gamma_1))^k \gamma_k.
   \label{e:p_symb_pf}
\end{equation}
\end{cor}

\section{The Case $G=\Sp(n,\R)$}

Next we consider the case where $G=\Sp(n,\R)$, or $\g=\sp_{n}$.
Let $E_{i,j}$ be as in the previous section.
Write an element $Z \in V=\u$ as
$Z=\sum_{i,j \in [n]} z_{i,j}E_{i,-j}$, $z_{j,i}=z_{i,j}$
and let
$$
\tilde{\pd}_{i,j}= \begin{cases}
                    2\pd_{i,j} &(i=j) \\ 
                    \pd_{i,j} &(i \ne j).
                  \end{cases}
$$
Let us denote by
$M,D$ symmetric matrices of size $n$
whose $(i,j)$-th entry are given by
the multiplication operators $z_{i,j}$ 
and the differential operators $\tilde{\pd}_{i,j}$.
Then
$$
\varGamma_k
  :=\sum_{ \substack{I, J \subset[n] \\ |I|=|J|=k} }  
    \det(z^{I}_{J})\det(\tilde{\pd}^{I}_{J})
  \quad (k=0,1,\dots,n)
$$
form a generating system for 
$\mathscr{PD}(V)^{\Kc}$ with $\Kc=\GL{n}$,
where
$z^{I}_{J},\pd^{I}_{J}$ denote submatrices of $M,D$
consisting of the entries
whose row- and column- indices are 
in $I$ and $J$, respectively
(\cite{HU91}).

Take a basis $\{ X_{i,j}^{\epsilon} \}_{\epsilon=0,\pm;i,j \in [n]}$
for $\sp_{n}$ as follows:
{\allowdisplaybreaks
\begin{equation}
 \begin{aligned}
\label{e:basis_sp}
 X^{0}_{i,j} &:= E_{i,j}-E_{-j,-i}  \quad ( 1 \leqsl i, j\leqsl n),
  \\
 X^{+}_{i,j} &:= E_{i,-j}+E_{j,-i} \quad (1 \leqsl i \leqsl j\leqsl n),
  \\
 X^{-}_{i,j} &:= E_{-j,i}+E_{-i,j} \quad (1 \leqsl i \leqsl j\leqsl n).
 \end{aligned}
\end{equation}
}

%
%
\begin{prop}
\label{p:rep_op_sp}
The differential operators $\d \pi_{\lambda}(X^{\epsilon}_{i,j})$
$(\epsilon=0,\pm;i,j \in [n])$
are given by
{\allowdisplaybreaks
\begin{align}
 \d \pi_{\lambda}(X^{0}_{i,j})
      &= s \delta_{i,j} -\sum_{k=1}^{n} z_{j,k} \tilde{\pd}_{i,k},
  \label{e:op_a_{i,j}_sp}    \\
 \d \pi_{\lambda}(X^{+}_{i,j})
      &= -\tilde{\pd}_{i,j}, 
  \label{e:op_b_{i,j}_sp}    \\
 \d \pi_{\lambda}(X^{-}_{i,j})
      &= -2s z_{i,j} -\sum_{1 \leqsl k,l \leqsl n} 
                        z_{k,i} z_{j,l} \tilde{\pd}_{k,l}.
  \label{e:op_c_{i,j}_sp}  
\end{align}
}
\end{prop}
\begin{proof}
Let $\umat{n}_{i,j}$ be as in
the proof of Proposition \ref{p:rep_op_so},
and let $\tilde{z} := \sum_{k,l} z_{k,l} \umat{n}_{k,l}$
with $z_{l,k}=z_{k,l}$.
Then $z:=\tilde{z} J_n$ belongs to $\Omega$
if it is positive definite.
Then, we can calculate the differential operators 
$\d \pi_{\lambda}(X^{\epsilon}_{i,j})$ 
as in the case of $\so_{2n}$.

(I) Writing 
$
\left[\begin{smallmatrix} a & 0 \\ 0 & d \end{smallmatrix}\right]
   := \exp (-t X^{0}_{ij}) 
$,
we see that
\begin{align*}
 \exp(-t X^{0}_{ij}).z  
  &= a z d^{-1}
     \\
  &= \sum_{k,l \in [n]} 
      \left( z_{k,l}-t(\delta_{i,k} z_{j,l} + \delta_{i,l} z_{k,j}) \right) 
            \umat{n}_{k,l} J_n + O(t^2) 
\end{align*}
and  
\[
 (\det d)^{s}
    =e^{s t \trace{\umat{n}_{n+1-j,n+1-i}} + O(t^2)},
\]
and hence obtain that
\begin{align*}
 \d \pi_{\lambda}(X^{0}_{i,j})
  &=  s \delta_{i,j}
       - \sum_{i \leqsl l} z_{j,l} \pd_{i,l} 
       - \sum_{k \leqsl i} z_{k,j} \pd_{k,i} 
     \\
  &=  s \delta_{i,j} -\sum_{k \in [n]} z_{j,k} \tilde{\pd}_{i,k}.
\end{align*}

(II) Writing 
$
\left[\begin{smallmatrix} 1 & b \\ 0 & 1 \end{smallmatrix}\right]
   := \exp (-t X^{+}_{ij})
$,
we see that
\begin{align*}
 \exp(-t X^{+}_{ij}).z 
  &= z + b 
     \\
  &= \left( \tilde z - t (\umat{n}_{i,j} + \umat{n}_{j,i}) \right) J_n
     \\
  &= \sum_{k,l \in [n]}
      \left(
        z_{k,l}-t \delta_{i,k} \delta_{j,l} - t \delta_{j,k} \delta_{i,l}
      \right) \umat{n}_{k,l} J_n,
\end{align*}
and hence obtain that
\begin{align*}
 \d \pi_{\lambda}(X^{+}_{i,j})
  &= - \sum_{k \leqsl l} 
       ( \delta_{i,k} \delta_{j,l} +\delta_{i,l} \delta_{j,k}) \pd_{k,l}
     \\
  &= - \tilde{\pd}_{i,j}.
\end{align*}

(III) Writing 
$
\left[\begin{smallmatrix} 1 & 0 \\ c & 1 \end{smallmatrix}\right]
   := \exp (-t X^{-}_{ij})
$,
we see that
\begin{align*}
 \exp(-t X^{-}_{ij}).z
  &=z (1+cz)^{-1}
     \\
  &= z \left( 1_n - t J_n (\umat{n}_{i,j}+\umat{n}_{j,i}) z \right)^{-1}
     \\
  &= \left( 
      \tilde{z} +t \tilde{z}(\umat{n}_{i,j}+\umat{n}_{j,i}) \tilde{z}
     \right) J_n + O(t^2)
     \\
  &= \sum_{k,l} 
       \left( z_{k,l} 
             + t ( z_{k,i} z_{j,l}+ z_{k,j} z_{i,l} )
       \right) \umat{n}_{k,l} J_n + O(t^2),
\end{align*}
and 
\begin{align*}
 \det(1+cz)^s 
  &=\left( 
      \det( 1_n - t J_n ( \umat{n}_{i,j}+\umat{n}_{j,i}) \tilde{z} J_n )
    \right)^s
  = e^{-st \trace{( \umat{n}_{i,j}+\umat{n}_{j,i}) \tilde{z}} +O(t^2) }
     \\
  &= e^{ -s t (z_{j,i}+z_{i,j}) + O(t^2) }
  =e^{-2 s t z_{i,j} + O(t^2) },
\end{align*}
and hence obtain that
\begin{align*}
 \d \pi_{\lambda}(X^{-}_{i,j}) 
  &=-2s z_{i,j} 
    - \sum_{k \leqsl l} 
          \left( z_{k,i}z_{j,l} + z_{k,j}z_{i,l} \right) \pd_{k,l}
     \\
  &=-2s z_{i,j} 
    - \sum_{k,l \in [n]} z_{k,i}z_{j,l} \tilde{\pd}_{k,l}.
\end{align*}
This completes the proof.
\end{proof}

Noting that the dual basis of \eqref{e:basis_sp}
is given by
\begin{align*}
 (X^{0}_{i,j})^{\vee} 
   &= X^{0}_{j,i},
       &   
 (X^{\pm}_{i,j})^{\vee}
   &= \begin{cases}
       X^{\mp}_{i,j}         & ( i \ne j),
            \\
       \frac12 X^{\mp}_{i,j} & ( i = j),
      \end{cases}
\end{align*}
let us define $\sigma(\bs{X})$ 
by substituting $\tilde{\xi}_{i,j}$ into $\tilde{\pd}_{i,j}$
following the prescription 
$\eqref{e:symbol_of_ex_Casimir}$
as above.

%
%
\begin{thm}
\label{t:nilpotency_sp}
Let
$u(z):=\exp \sum_{i \leqsl j} z_{i,j} X^{+}_{i,j} \in U_{\Omega}$.
Then we have
\begin{align}
 & \Ad(u(z)^{-1}) \sigma(\bs{X})
    = s \sum_{i} X^{0}_{i,i}
       - \sum_{i \leqsl j} \tilde\xi_{i,j} X^{-}_{i,j}
     \\
 &=\left[
   \begin{array}{cccc|cccc}
   s &               &           &     
       &\hph{\xi_{1,n-1}} &\hph{\xi_{2,n-1}} &\hph{\cdots}&\hph{\xi_{2,n-1}}
    \\
     &        s     &            & 
        &               &           &            & 
    \\ 
     &              & \ddots     & 
        &               &           &            & 
    \\ 
     &              &            & s
        &               &           &             & 
    \\[8pt]
   \hline
   -\xi_{1,n-1}  & -\xi_{2,n-1}   & \cdots & -2\xi_{n,n}    
        & -s            &           &             &
    \\
   \vdots        & \vdots    & \adots        & \vdots 
        &               & \ddots    &             &
    \\
   -\xi_{1,2}    & -2\xi_{2,2} & \cdots     & -\xi_{2,n-1}
        &           &           & -s          &
    \\
   -2\xi_{1,2}   & -\xi_{1,2} & \cdots     & -\xi_{1,n-1}    
        &          &           &             & -s
   \end{array}
   \right].
\end{align} 
\end{thm}
\begin{proof}
Again, this is just a simple matrix calculation,
and can be shown in the same way 
as in the case of $\so_{2n}$.
In fact, if we write
$
\left[
   \begin{smallmatrix}
     A & B \\
     C & D   
   \end{smallmatrix}
\right]
 := \sigma_{\lambda}(\bs{X})
$,
then 
the summation 
$\sum_{k} z_{i,k} \tilde{\xi}_{j,k}$ 
in $\sigma_{\lambda}(X^{0}_{j,i})=A_{ij}$
equals  $-(zC)_{ij}$, 
from which it follows that
\[
 A = s 1_n + zC.
\]
Similarly,
the summation 
$\sum_{k,l} z_{i,k} z_{l,j} \tilde{\xi}_{k,l}$ 
in $\sigma_{\lambda}(X^{-}_{i,j})=(B J_n)_{ij}$
equals $(zCz J_n)_{ij}$,
from which it follows that
\[
 B = -2s z - zCz.
\]
Now exactly the same calculation 
as in the proof of Theorem \ref{t:nilpotency_so}
yields the identity to be shown.
\end{proof}

As in the previous section,
we rewrite 
$\sigma(X^{0}_{i,i})$ and $\sigma(X^{-}_{i,j})$
using 
$\gamma_1=\sum_{k,l \in [n]} z_{k,l} \tilde\xi_{k,l}$ 
and substitute $s-u$ into $\gamma_1$ 
with a new indeterminate $u$,
which we denote by
$\tild{\sigma}(X^{0}_{i,i})$ and $\tild{\sigma}(X^{-}_{i,j})$,
respectively:
\begin{align}
 \tild{\sigma}(X^{0}_{i,i}) 
  &:= u 
      +\sum_{ \substack{ k,l \in [n] \\ l \ne i} } 
                  z_{k,l} \tilde{\xi}_{k,l},
   \label{tilde_sigma_X^{0}_{ii}_sp}
     \\
 \tild{\sigma}(X^{-}_{i,j})
  &:= -(u+s) z_{i,j} 
     - \sum_{k,l \in [n]}
        (z_{k,i}z_{j,l} - z_{i,j} z_{k,l} ) \tilde{\xi}_{k,l}.
   \label{tilde_sigma_X^{-}_{ij}_sp}
\end{align}
As for the others,
set $\tild{\sigma}(X^{\epsilon}_{i,j}):=\sigma(X^{\epsilon}_{i,j})$.
Then we define $\tild{\sigma}(\bs{X})$ by
\[
\tild{\sigma}(\bs{X})
 :=\sum_{\epsilon;i,j} 
      \tild{\sigma}((X^{\epsilon}_{i,j})^{\vee}) \otimes X^{\epsilon}_{i,j}.
\]

%
%
\begin{thm}
\label{t:moment_map_with_u_sp}
Let $u(z) \in U_{\Omega}$ be as in Theorem \ref{t:nilpotency_sp},
then we have
\begin{align}
 &\Ad(u(z)^{-1}) \tild{\sigma}(\bs{X})
  =(u+\gamma_1) \sum_{i} X^{0}_{i,i}
    - \sum_{i \leqsl j} \tilde\xi_{i,j} X^{-}_{i,j}
    - (s-(u+\gamma_1)) \sum_{i \leqsl j} z_{i,j} X^{+}_{i,j}
      \\
 &=\left[
   \begin{array}{cccc|cccc}
   u+\gamma_1 &               &           &     
       &-\tau\, z_{1,n-1}       & \cdots    &-\tau\, z_{1,2}    &-\tau\,z_{1,1}
    \\
     & u+\gamma_1  &            & 
        &-\tau\,z_{2,n-1}      & \cdots    &-\tau\,z_{2,2}      &-\tau\,z_{1,2}
    \\ 
     &              & \ddots     & 
        & \vdots        & \adots    & \vdots     & \vdots
    \\ 
     &              &            & u+\gamma_1
        &-\tau\,z_{1,2} & \cdots    &-\tau\,z_{2,n-2} &-\tau\,z_{1,n-1}
    \\
   \hline
   -\xi_{1,n-1}  & -\xi_{2,n-1}   & \cdots & -2\xi_{n,n}    
        & -u-\gamma_1   &           &             &
    \\
   \vdots        & \vdots    & \adots        & \vdots 
        &               & \ddots    &             &
    \\
   -\xi_{1,2}    & -2\xi_{2,2} & \cdots     & -\xi_{2,n-1}
        &           &           &-u-\gamma_1  &
    \\
   -2\xi_{1,1}   & -\xi_{1,2} & \cdots     & -\xi_{1,n-1}    
        &          &           &             & -u-\gamma_1
   \end{array}
   \right].
      \label{e:sigX2_sp}
\end{align} 
Here, we set $\tau:=s-(u+\gamma_1)$ 
in \eqref{e:sigX2_sp} for brevity.
\end{thm}
\begin{proof}
The theorem follows from matrix calculation 
parallel to that in the proof of Theorem \ref{t:moment_map_with_u}.

In fact,
if we write
$
\left[ 
   \begin{smallmatrix} A(u) & B(u) \\ C & D(u) \end{smallmatrix}
\right]
 := \tild\sigma(\bs{X}),
$
then we see that 
\begin{equation}
\label{e:A(u)_{ij}_sp}
 A(u)_{ij}
 = (u + \sum_{k,l} z_{k,l} \tilde{\xi}_{k,l} ) \delta_{i,j}
    - \sum_{k=1}^{n} z_{k,i} \tilde{\xi}_{k,j}
\end{equation}
and the second summation in \eqref{e:A(u)_{ij}_sp} 
equals 
$-(\tilde{z}\tilde{C})_{ij}=-(zC)_{ij}$,
hence obtain that 
\[
 A(u) = (u +\gamma_1) 1_n + zC.
\]
Similarly,
we see that the summation in 
$\tild{\sigma}(X^{-}_{i,j})=\tilde{B}(u)_{ij}$
equals $z_{i,j} \gamma_1+(\tilde{z} \tilde{C} \tilde{z})_{ij}$
as shown in the proof of Theorem \ref{t:nilpotency_sp},
hence obtain that
\begin{equation}
\label{e:keyformula_(2,2)-block_in_sigma_with_u_sp}
 B(u) = -(u+s+\gamma_1) z - zCz.
\end{equation}
Again, exactly the same matrix calculation 
as in the proof of Theorem \ref{t:moment_map_with_u}
yields the identity to be shown.
\end{proof}
It follows immediately from Theorem \ref{t:moment_map_with_u_sp}
and the formula \eqref{e:msf_det}
that determinant of $\tild\sigma(\bs{X})$ yields 
a generating function for $\{ \gamma_k \}$.
%
\begin{cor}
\label{c:gen_fun_sp}
Retain the notation above.
Then we have the following formula:
\begin{equation}
\label{e:p_symb_symm_capelli}
  \det(\tild\sigma(\bs{X}))
=(-1)^n \sum_{k=0}^{n}
   (u+\gamma_1)^{2n-2k} (s-(u+\gamma_1))^k \gamma_k.
\end{equation}
\end{cor}

\section{The Case $G=\SU(p,q)$ with $p \geqsl q$}

Finally, we consider the case where $G=\SU(p,q)$, 
or $\g=\sl_{p+q}$, with $p \geqsl q$.
In this section,
let $E_{i,j}$ denote the $(p+q) \times (p+q)$ matrix 
with its $(i,j)$-th entry being $1$ 
and all the others $0$.
Write an element $V=\u$ as
$Z=\sum_{i \in [p],j \in [q]} z_{i,j}E_{i,j}$.
Let us denote by
$M,D$ $p \times q$ matrices 
whose $(i,j)$-th entries are given by
the multiplication operators $z_{i,j}$ 
and the differential operators ${\pd}_{i,j}$.
Then
$$
\varGamma_k
  :=\sum_{ \substack{I \subset[p], J \subset[q] \\ |I|=|J|=k} } 
    \det(z^{I}_{J})\det({\pd}^{I}_{J})
  \quad (k=0,1,\dots,q)
$$
form a generating system for 
$\mathscr{PD}(V)^{\tilde{K}_{\C}}$ 
with 
$\tilde{K}_{\C}=\GL{p}\times \GL{q} \simeq \Kc \times \C^{\times}$,
where
$z^{I}_{J},\pd^{I}_{J}$ denote submatrices of $M,D$
consisting of the entries
whose row- and column- indices are 
in $I$ and $J$, respectively 
(\cite{HU91}).

Take a basis $\{ H_{i}, E^{\pm}_{i,j}, X^{\pm}_{i,j} \}$
for $\sl_{p+q}$ as follows:
{\allowdisplaybreaks
\begin{equation}
 \begin{aligned}
\label{e:basis_sl}
 H_{i} &:= E_{i,i}-E_{p+q,p+q} & &
        ( 1 \leqsl i \leqsl p+q-1) 
   \\
 E^{+}_{i,j} &:= E_{i,j}     &  &
        ( 1 \leqsl i \ne j \leqsl p)
   \\
 E^{-}_{i,j} &:= E_{p+i,p+j} &  &
        ( 1 \leqsl i \ne j \leqsl q)
   \\
 X^{+}_{i,j} &:= E_{i,p+j}   &  &
        (1 \leqsl i \leqsl p, 1 \leqsl j \leqsl q)
   \\
 X^{-}_{i,j} &:= E_{p+j,i}   &  &
        (1 \leqsl i \leqsl p, 1 \leqsl j \leqsl q)
 \end{aligned}
\end{equation}
}

%
%
\begin{prop}
The differential operators 
$
\d \pi_{\lambda}(H_{i}), %
\d \pi_{\lambda}(E^{\pm}_{i,j}), %
\d \pi_{\lambda}(X^{\pm}_{i,j})
$
are given by
{\allowdisplaybreaks
\begin{align}
 \d \pi_{\lambda}(H_{i})
    &= s -\sum\limits_{l=1}^{q} z_{i,l} {\pd}_{i,l} 
          -\sum\limits_{k=1}^{p} z_{k,q} \pd_{k,q}   &
                 & (1 \leqsl i \leqsl p)   
  \label{e:op_h_{i}_sl}                
            \\   
 \d \pi_{\lambda}(H_{p+j})
    &= \sum\limits_{k=1}^{p} 
            (z_{k,j} {\pd}_{k,j} - z_{k,q} \pd_{k,q}) &
                 & (1 \leqsl j \leqsl q-1)
  \label{e:op_h_{p+j}_sl}                
            \\
 \d \pi_{\lambda}(E^{+}_{i,j})
    &= -\sum_{l=1}^{q} z_{j,l} {\pd}_{i,l}  & 
        & ( 1 \leqsl i \ne j \leqsl p)
  \label{e:op_e^{+}_{i,j}_sl}    
            \\
 \d \pi_{\lambda}(E^{-}_{i,j})
    &= \sum_{k=1}^{p} z_{k,i}  {\pd}_{k,j}  & 
        & ( 1 \leqsl i \ne j \leqsl q)
  \label{e:op_e^{-}_{i,j}_sl}    
            \\
 \d \pi_{\lambda}(X^{+}_{i,j}) 
    &= -{\pd}_{i,j}                         & 
        & (1 \leqsl i \leqsl p, 1 \leqsl j \leqsl q)
  \label{e:op_x^{+}_{i,j}_sl}    
            \\
 \d \pi_{\lambda}(X^{-}_{i,j})
      &= -s z_{i,j} + \sum_{\substack{1 \leqsl k \leqsl p, \\%
                                      1 \leqsl l \leqsl q} 
                        } 
                         z_{k,j} z_{i,l} {\pd}_{k,l}.  & 
        & (1 \leqsl i \leqsl p, 1 \leqsl j \leqsl q)
  \label{e:op_x^{-}_{i,j}_sl}  
\end{align}
} 
\end{prop}
\begin{proof}
Let us denote the $p \times q$ matrix 
whose $(i,j)$-th entry is $1$ and all the others are $0$
by $\umat{p,q}_{i,j}$,
and let 
$z:=\sum_{i \in [p], j \in [q]} z_{i,j} \umat{p,q}_{i,j}$.
Then $z$ belongs to $\Omega$ if it is positive definite.

(I) 
Write 
$
\left[\begin{smallmatrix} a & 0 \\ 0 & d \end{smallmatrix}\right]
:=\exp (-t H_{i}) 
$.
For $i=1,\dots,p$, 
we see that
\begin{align*}
 \exp(-t H_{i}).z  
  &= a z d^{-1}
     \\
  &= \left(1_p -t \umat{p}_{i,i}+O(t^2)\right) z
       \left(1_q+t \umat{q}_{q,q} +O(t^2) \right)^{-1}
     \\
  &= \sum_{k \in [p], l \in [q]} 
      \left( z_{k,l}-t(\delta_{i,k} z_{i,l} + \delta_{l,q} z_{k,q}) \right) 
            \umat{p,q}_{k,l} + O(t^2) 
\end{align*}
and  
\[
 (\det d)^{s}
    =e^{s t \trace{\umat{q}_{p+q,p+q}} + O(t^2)},
\]
and hence obtain that
\begin{equation*}
 \d \pi_{\lambda}(H_{i})
  =  s 
       - \sum_{l \in [q]} z_{i,l} \pd_{i,l} 
       - \sum_{k \in [p]} z_{k,q} \pd_{k,q} 
   \quad (i=1,\dots,p).
\end{equation*}
For $j=1,\dots,q-1$, 
we see that
\begin{align*}
 \exp(-t H_{p+j}).z  
  &= a z d^{-1}
     \\
  &= z \left(1_q -t (\umat{q}_{j,j}-\umat{q}_{q,q}) +O(t^2) \right)^{-1}
     \\
  &= \sum_{k \in [p], l \in [q]} 
      \left( z_{k,l} + t (\delta_{l,j} z_{k,j} - \delta_{l,q} z_{k,q}) \right) 
            \umat{p,q}_{k,l} + O(t^2) 
\end{align*}
and $(\det d)^{s}=1$,
hence obtain that
\begin{equation*}
 \d \pi_{\lambda}(H_{p+j})
  = \sum_{k \in [p]} (z_{k,j} \pd_{k,j}- z_{k,q} \pd_{k,q})
     \quad (j=1,\dots,q-1). 
\end{equation*}

(II) Writing 
$
\left[\begin{smallmatrix} a & 0 \\ 0 & 1 \end{smallmatrix}\right]
:= \exp (-t E^{+}_{ij}) 
$,
we see that
\begin{align*}
 \exp(-t E^{+}_{i,j}).z  
  &= a z
     \\
  &= \left( 1_p -t \umat{p}_{i,j} \right) z
     \\
  &= \sum_{k \in [p], l \in [q]} 
      \left( z_{k,l}-t \delta_{i,k} z_{j,l} \right) \umat{p,q}_{k,l} + O(t^2) 
\end{align*}
and hence obtain that
\begin{equation*}
 \d \pi_{\lambda}(E^{+}_{i,j})
  = -\sum_{l \in [q]} z_{j,l} \pd_{i,l}.
\end{equation*}

(III) Writing 
$
\left[\begin{smallmatrix} 1 & 0 \\ 0 & d \end{smallmatrix}\right]
:= \exp (-t E^{-}_{ij}) 
$,
we see that
\begin{align*}
 \exp(-t E^{-}_{i,j}).z  
  &= z d^{-1}
     \\
  &= z \left( 1_q -t \umat{q}_{i,j} \right)^{-1}
     \\
  &= \sum_{k \in [p], l \in [q]} 
      \left( z_{k,l} + t \delta_{j,l} z_{k,i} \right)\umat{p,q}_{k,l} + O(t^2) 
\end{align*}
and $\det d=1$,
and hence obtain that
\begin{equation*}
 \d \pi_{\lambda}(E^{-}_{i,j})
  = \sum_{k \in [p]} z_{k,i} \pd_{k,j}.
\end{equation*}

(IV) Writing 
$
\left[\begin{smallmatrix} 1 & b \\ 0 & 1 \end{smallmatrix}\right]
:=\exp (-t X^{+}_{ij})
$,
we see that
\begin{align*}
 \exp(-t X^{+}_{ij}).z 
  &= z + b 
     \\
  &= z - t \umat{p,q}_{i,j}
\end{align*}
and hence obtain that
\begin{equation*}
 \d \pi_{\lambda}(X^{+}_{i,j})
  = - \pd_{i,j}.
\end{equation*}

(V) Writing 
$
\left[\begin{smallmatrix} 1 & 0 \\ c & 1 \end{smallmatrix}\right]
 := \exp (-t X^{-}_{ij})
$,
we see that
\begin{align*}
 \exp(-t X^{-}_{ij}).z
  &=z (1+cz)^{-1}
     \\
  &= z \left( 1_q - t \umat{q,p}_{j,i} z \right)^{-1}
     \\
  &= z + t z \umat{q,p}_{j,i} z + O(t^2)
     \\
  &= \sum_{k \in [p], l \in [q]} 
       \left( 
          z_{k,l} + t  z_{k,j} z_{i,l}
       \right) \umat{p,q}_{k,l} + O(t^2),
\end{align*}
and 
\begin{align*}
 \det(1+cz)^s 
  &=\left( 
      \det( 1_q - t \umat{q,p}_{j,i} z )
    \right)^s
  = e^{-st \trace{\umat{q,p}_{j,i} z} +O(t^2) }
     \\
  &= e^{ -s t z_{i,j} + O(t^2) },
\end{align*}
and hence obtain that
\begin{equation*}
 \d \pi_{\lambda}(X^{-}_{i,j}) 
  =-s z_{i,j} 
    - \sum_{k \in [p], l \in [q]} 
          z_{k,j} z_{i,l} \pd_{k,l}.
\end{equation*}
This completes the proof.
\end{proof}

Noting that the dual basis of \eqref{e:basis_sl}
is given by
\begin{align}
 H_{i}^{\vee} 
   &= H_{i}-\tfrac{1}{p+q} \sum_{k=1}^{p+q-1} H_{k}, 
     &  
 (E^{\pm}_{i,j})^{\vee} 
   &= E^{\pm}_{j,i},
     &  
 (X^{\pm}_{i,j})^{\vee} 
   &= X^{\mp}_{i,j},
     \notag 
\end{align}
let us define $\sigma(\bs{X})$ 
following the prescription 
$\eqref{e:symbol_of_ex_Casimir}$ 
as above.

%
%
\begin{thm}
\label{t:nilpotency_sl}
Let
$
u(z):=\exp \sum_{i \in [p], j \in [q]} z_{i,j} X^{+}_{i,j} %
 \in U_{\Omega}
$
Then we have
\begin{align}
 & \Ad(u(z)^{-1}) \sigma(\bs{X})
    = \tfrac{p}{p+q}s \sum_{i=1}^{p} H_{i} 
       -\tfrac{q}{p+q}s \sum_{j=1}^{q-1} H_{p+j}
       - \sum_{i \leqsl j} \xi_{i,j} X^{-}_{i,j}
     \\
 &=\left[
   \begin{array}{cccc|cccc}
  \frac{p}{p+q}s &               &           &     
       &\hph{-\xi_{1,q}} &\hph{-\xi_{2,q}} &\hph{\cdots}&\hph{-\xi_{p,q}}
    \\
     & \frac{p}{p+q}s     &            & 
        &               &           &            & 
    \\ 
     &              & \ddots     & 
        &               &           &            & 
    \\ 
     &              &            &  \frac{p}{p+q}s
        &               &           &             & 
    \\[8pt]
   \hline
   -\xi_{1,1}  & -\xi_{2,1}   & \cdots   & -\xi_{p,1}    
        & -\frac{q}{p+q}s &      &          &
    \\
   -\xi_{1,2}  & -\xi_{2,2}   & \cdots   & -\xi_{p,2}
        &           & -\frac{q}{p+q}s&      &
    \\
   \vdots        & \vdots    &           & \vdots 
        &               &         & \ddots      &
    \\
   -\xi_{1,q}   & -\xi_{2,q} & \cdots    & -\xi_{p,q}    
        &          &           &             & -\frac{q}{p+q}s
   \end{array}
   \right].
\end{align} 
\end{thm}
\begin{proof}
It follows from \eqref{e:op_h_{i}_sl} and 
\eqref{e:op_h_{p+j}_sl} that
\[
 \frac{1}{p+q} \sum_{k=1}^{p+q-1} \sigma(H_{k})
  = \frac{p}{p+q}s  - \sum_{k \in [p]} z_{k,q} \xi_{k,q}.
\]
Thus we obtain that
\begin{equation}
\label{e:symbol_h_{i}_sl}
 \begin{aligned}
 \sigma(H_{i}^{\vee})
   & = \frac{q}{p+q} s - \sum_{l \in [q]} z_{i,l} \xi_{i,l}  &
        & (i=1,\dots,p)
          \\
 \sigma(H_{p+j}^{\vee})
   & = -\frac{p}{p+q} s + \sum_{k \in [p]} z_{k,j} \xi_{k,j}    &
        & (j=1,\dots,q-1).
\end{aligned}
\end{equation}
Now, if we write
$
\left[
   \begin{smallmatrix}
     A & B \\
     C & D   
   \end{smallmatrix}
\right]
 := \sigma_{\lambda}(\bs{X})
$,
then 
\begin{subequations} 
\label{e:sigX_sl}
 \begin{align}
 A_{ij} 
  &= \frac{q}{p+q}s \delta_{i,j} -\sum_{l \in [q]} z_{i,l} \xi_{j,l}
     &  &(i,j \in [p]) 
     \label{e:A_{ij}_sl}
     \\
 B_{ij} 
  &= -s z_{i,j} 
     + \sum_{k \in [p], j \in [q]} z_{k,j}z_{i,l}\xi_{k,l}
     &  &(i \in [p],j \in [q]) 
     \label{e:B_{ij}_sl}
     \\
 C_{ij}
  &= -\xi_{j,i}
     &  &(i \in [p], j \in [q])
     \label{e:C_{ij}_sl}
     \\
 D_{ij}
  &=-\frac{p}{p+q}s \delta_{i,j} +\sum_{k \in [p]} z_{k,j} \xi_{k,i} 
     &  &(i,j \in [q])
     \label{e:D_{ij}_sl}
 \end{align}
\end{subequations}
by \eqref{e:op_e^{+}_{i,j}_sl}, 
\eqref{e:op_e^{-}_{i,j}_sl}, \eqref{e:op_x^{+}_{i,j}_sl},
\eqref{e:op_x^{-}_{i,j}_sl} and \eqref{e:symbol_h_{i}_sl}.
Now exactly the same matrix calculation 
as in the proof of Theorem \ref{t:nilpotency_sp}
yields the formula to be shown.
\end{proof}

As in the previous sections,
we rewrite 
$\sigma(H_{i}^{\vee})$ and $\sigma(X^{-}_{i,j})$ 
using 
$\gamma_1=\sum_{k \in [p], l \in [q]} z_{k,l} \xi_{k,l}$ 
and substitute $s-u$ into $\gamma_1$ 
with a new indeterminate $u$,
which we denote by
$\tild{\sigma}(H_{i}^{\vee})$ 
and $\tild{\sigma}(X^{-}_{i,j})$,
respectively:
\begin{align}
 \tild\sigma(H_{i}^{\vee})
  &:= (u -\frac{p}{p+q} s) 
     + \sum_{\substack{k \in [p], l \in [q] \\ k \ne i}} z_{k,l} \xi_{k,l}   &
  &(i=1,\dots,p )
        \\
 \tild\sigma(H_{p+j}^{\vee})
  &:= (-u +\frac{q}{p+q} s) 
     + \sum_{\substack{k \in [p], l \in [q] \\ l \ne j}} z_{k,l} \xi_{k,l}   &
  &(j=1,\dots, q-1)
     \\
 \tild\sigma(X^{-}_{i,j})
   &:= -u z_{i,j} 
       -\sum_{\substack{k \in [p],l \in [q] \\ k \ne i,l \ne j}} 
             (z_{i,j}z_{k,l}-z_{i,l}z_{k,j}) \xi_{k,l}.     &
   &(i=1,\dots,p; j=1,\dots,q)
\end{align}
As for the others,
set $\tild{\sigma}(\cdot):=\sigma(\cdot)$.
Then we define $\tild{\sigma}(\bs{X})$ by
\[
\tild{\sigma}(\bs{X})
 =\sum_{i} \tild{\sigma}(H^{\vee}_{i}) \otimes H_{i}
  + \sum_{\epsilon;i,j} 
      \left(
       \tild{\sigma}((E^{\epsilon}_{i,j})^{\vee}) \otimes E^{\epsilon}_{i,j}
     + \tild{\sigma}((X^{\epsilon}_{i,j})^{\vee}) \otimes X^{\epsilon}_{i,j}
      \right).
\]

%
%
\begin{thm}
\label{t:moment_map_with_u_sl}
Let $u(z) \in U_{\Omega}$ be as in Theorem \ref{t:nilpotency_sl},
then we have
\begin{align}
 &\Ad(u(z)^{-1}) \tild{\sigma}(\bs{X})
  =\left( u+\gamma_1 -\tfrac{p}{p+q}s \right) \sum_{i=1}^{p} H_{i}
   +\left(-u-\gamma_1 +\tfrac{q}{p+q}s \right) \sum_{j=1}^{q-1} H_{p+j}
      \notag     \\
 & \hph{\Ad(u(z)^{-1}) \tild{\sigma}(\bs{X})=}\;
    - \sum_{i,j} \xi_{i,j} X^{-}_{i,j}
    - (s-(u+\gamma_1)) \sum_{i,j} z_{i,j} X^{+}_{i,j}
                 \\
 &=\left[
   \begin{array}{cccc|cccc}
  u_{+} &              &           &     
       &-\tau\,z_{1,1}        &-\tau\,z_{2,1} & \cdots     &-\tau\,z_{1,q}
    \\
               &u_{+}  &           & 
        &-\tau\,z_{2,1}       &-\tau\,z_{2,2} & \cdots     &-\tau\,z_{2,q}
    \\ 
     &              & \ddots     & 
        & \vdots        & \vdots    &         & \vdots
    \\ 
     &              &            & u_{+}
        &-\tau\,z_{p,1}       &-\tau\,z_{p,2} & \cdots     &-\tau\,z_{p,q}
    \\
   \hline
   -\xi_{1,1}    & -\xi_{2,1}      & \cdots       & -\xi_{p,1}    
        & u_{-}     &                 &             &
    \\
   -\xi_{1,2}    & -\xi_{2,2}      & \cdots       & -\xi_{p,2}
        &           & u_{-}           &             &
    \\
   \vdots        & \vdots    &     & \vdots 
        &           &                 & \ddots      &
    \\
   -\xi_{1,q}   & -\xi_{2,q}       & \cdots       & -\xi_{p,q}    
        &          &           &             & u_{-}
   \end{array}
   \right].
      \label{e:sigX2_sl}
\end{align} 
Here, we set 
$u_{+}:=u+\gamma_1 -\tfrac{p}{p+q}s, %
u_{-}:=-u-\gamma_1 +\tfrac{q}{p+q}s$
and $\tau:=s-(u+\gamma_1)$
in \eqref{e:sigX2_sl} for brevity.
\end{thm}
\begin{proof}
If we write
$
\left[ 
  \begin{smallmatrix} A(u) & B(u) \\ C & D(u) \end{smallmatrix}
\right]
 := \tild\sigma(\bs{X})
$,
then we can show that
\begin{align}
 A(u)
  &= (u +\gamma_1 -\tfrac{p}{p+q} s) 1_p + zC
     \\
 B(u)
  &=-(u+\gamma_1) z-zCz 
     \\
 D(u)
  &= (-u -\gamma_1 +\tfrac{q}{p+q} s) 1_q -Cz
 \end{align}
exactly in the same way as 
in the proof of Theorems \ref{t:moment_map_with_u}
and \ref{t:moment_map_with_u_sp},
and hence obtain the theorem.
\end{proof}
It follows immediately from Theorem \ref{t:moment_map_with_u_sl}
and the formula \eqref{e:msf_det}
that determinant of $\tild\sigma(\bs{X})$ yields 
a generating function for $\{ \gamma_k \}$.
%
%
\begin{cor}
\label{c:gen_fun_sl}
Retain the notation above.
Then we have the following formula:
\begin{equation}
\label{e:p_symb_capelli}
  \det(\tild\sigma(\bs{X}))
= (-1)^q \sum_{k=0}^{q}
            (u+\gamma_1-\tfrac{p}{p+q}s)^{p-k}%
            (u+\gamma_1-\tfrac{q}{p+q}s)^{q-k}%
            (s-(u+\gamma_1))^k \gamma_k.
\end{equation}
\end{cor}

%
%

\section{Concluding Remark}

In general,
let $G$ be a Lie group,
$\g$ its Lie algebra, and $\g^{*}$ the dual of $\g$. 
If $M$ is a $G$-manifold,
then the cotangent bundle $T^{*}M$ 
is a symplectic $G$-manifold,
and the moment map $\mu:T^{*}M \to \g^{*}$
is given by
\begin{equation}
\label{e:def_of_moment_map}
 \langle \mu(x,\xi),X \rangle = \xi(X_{M}(x))
         \qquad  
      (x \in M, \xi \in T^{*}_{x} M) 
\end{equation}
for $X \in \g$,
where $\langle \cdot,\cdot \rangle$ denotes 
the canonical pairing between $\g^{*}$ and $\g$,
and $X_{M}(x) \in T_x M$ the tangent vector at $x \in M$
defined by
\begin{equation}
\label{e:def_vec_field}
 X_{M}(x) f := \left. \frac{d}{dt} \right|_{t=0} f( \exp(-t X).x ) 
\end{equation}
for functions $f$ defined around $x \in M$ 
(see e.g.~\cite{MFK94}).

Returning to our case,
if the character $\lambda$ is trivial,
it follows from \eqref{e:def_of_moment_map} and 
\eqref{e:def_vec_field} that
$\sigma(\bs{X})$ is identical to 
the moment map $\mu: T^{*}(\Gc/P) \to \g^{*}$
composed by the isomorphism $\g^{*} \simeq \g$ 
via Killing form
(or the $\Gc$-invariant nondegenerate bilinear form 
$B$ given in \eqref{e:bilinear_form}; 
our Lie algebras are simple!), 
which we also denote by $\mu$.

If $\lambda$ is nontrivial,
then $\sigma(\bs{X})$ can be regarded as 
a variant of the twisted moment map 
$\mu_{\lambda} : T^{*}(\Gc/P) \to \g^{*} \simeq \g$ 
(\cite{SV96});
the difference $\mu_{\lambda}-\mu$,
which they denote by $\lambda_{x}$ with $x \in \Gc/P$ therein,
can be expressed as $\lambda_{x}=\Ad(g) \lambda^{\vee}$,
where 
$\lambda^{\vee} \in \g$ corresponds to $\lambda \in \g^{*}$ 
via the nondegenerate bilinear form $B$
and 
$g$ is an element of a compact real form $U_{\R}$ of $\Gc$
such that $x=g \dot{e}$
with $\dot{e}$ the origin $P$ of $\Gc/P$.
Note that if $x$ is in the open subset 
$G/K=GP/P=U_{\Omega}P/P \subset \Gc/P$,
one can choose $u_x$ from $U_{\Omega}$
so that $x=u_x \dot{e}$
instead of $g$ from $U_{\R}$,
which was denoted by $u(z)$ in \S\S 4-6.
Then it is immediate to verify that
\begin{equation}
\sigma(\bs{X})=\Ad(u_x)\lambda^{\vee}+\mu(x,\xi).
   \notag 
\end{equation}
Writing  $\mu'_{\lambda}(x,\xi):=\sigma(\bs{X})$
to make its dependence on $x$ and $\xi$ transparent,
Theorems \ref{t:nilpotency_so}, \ref{t:nilpotency_sp},
and \ref{t:nilpotency_sl} state that an analogue of 
$U_{\R}$-equivariance of the twisted moment map holds:
\begin{equation}
\label{e:equivariance}
 \Ad(u_x^{-1})\mu'_{\lambda} (x,\xi)
 = \mu'_{\lambda} (u_x^{-1}.(x,\xi))
\end{equation}
with $x=u_x \dot{e}$,
though $U_{\Omega}$ is not a subgroup.
Since $u_x^{-1}.(x,\xi)=(\dot{e},\xi)$,
the relation \eqref{e:equivariance} says that
all the principal symbols of the representation operators 
of elements of the Lie algebras 
are determined by those at the origin.

At this moment, 
the following question naturally arise:
\textit{
What is the geometric meaning of the indeterminate $u$
introduced in the definition of $\tild{\sigma}(\bs{X})$?
}

We will seek for the answer to this question 
in the forthcoming paper.

%
%
\vspace{12pt}
\noindent
\textbf{Acknowledgments}\quad
I would like to thank T{\^o}ru Umeda 
for guiding me to the world of the Capelli identity, 
as well as for his valuable comments and advice.

\appendix

\section{Minor Summation Formulae}

In this appendix,
we collect formulae concerning Pfaffian and determinant
needed to show Corollaries \ref{c:gen_fun_so},
\ref{c:gen_fun_sp} and \ref{c:gen_fun_sl}.

Let $X$ be a $2n \times 2n$ matrix 
alternating along the anti-diagonal.
Since $X J_{2n}$ is an alternating matrix,
one can define Pfaffian of $X J_{2n}$,
which is denoted by $\Pf{X}$ 
as mentioned above.
If we write $X$ as 
$
X=\left[ \begin{smallmatrix}
          a & b \\ c & -J_n \tp{a} J_n
         \end{smallmatrix}
  \right]
$
with submatrices $a,b,c$ 
all being of size $n \times n$,
then $b$ and $c$ are also 
alternating along the anti-diagonal.
Thus we parametrize the submatrices 
$a,b,c$ as follows:
\begin{equation*}     
\label{parametrization} 
 a=\begin{bmatrix}
     a_{1,1} & \cdots & a_{1,n} \\
     \vdots  &        & \vdots \\
     a_{n,1} & \cdots & a_{n,n}
   \end{bmatrix},
\quad
b=\begin{bmatrix}
    b_{1,n} & \cdots & b_{1,2} & 0         \\
    \vdots  & \adots &   0     & -b_{1,2}  \\[-6pt]
            &        &         &           \\[-2pt]
    b_{n-1,n}&  0    & \adots  & \vdots    \\[6pt]
      0     &  -b_{n-1,n} & \cdots & -b_{1,n}
  \end{bmatrix},
\quad
c=\begin{bmatrix}
    c_{1,n} & \cdots & c_{n-1,n} & 0         \\
    \vdots  & \adots &   0     & -c_{n-1,n}  \\[-6pt]
            &        &         &           \\[-2pt]
    c_{1,2} &  0    & \adots  & \vdots    \\[6pt]
      0     &  -c_{1,2} & \cdots & -c_{1,n}
  \end{bmatrix}.
\end{equation*}
\begin{thm}
[\cite{IW06}]
\label{t:msf_pf_square}
Let 
$
X=\left[ \begin{smallmatrix}
          a & b \\ c & -J_n \tp{a} J_n
         \end{smallmatrix}
  \right]
$ 
be a $2n \times 2n$ matrix 
alternating along the anti-diagonal
with submatrices $a,b,c$ parametrized as above.
Then we can expand Pfaffian $\Pf{X}$ 
in the following way.
\begin{equation}
\label{e:msf_pf_square}
  \Pf{X} = \sum_{k=0}^{ \lfloor n/2 \rfloor } 
          \sum_{ \substack{
                 I, J \subset [n] \\%
                 |I|=|J|=2k}  }
           \sgn{\bar{I},I} \sgn{\bar{J},J}
           \det(a^{\bar I}_{\bar J}) \Pf{b_{I}} \Pf{c_{J}}.
\end{equation}
Here, for $I,J \subset [n]$,
$\bar I$ denotes the complement of $I$ in $[n]$, 
$a^{I}_{J}$ the submatrix of $a$ 
whose row- and column-indices
are in $I$ and $J$ respectively,
$b_{I}, c_{I}$ the submatrices of $b, c$
whose row- and column-indices are both in $I$,
and $\sgn{I,J}$ the signature of the permutation
$
\left(
  \begin{smallmatrix} 1,2,\cdots,n \\ I \;\; J \end{smallmatrix}
\right)
$.
\end{thm}

In order to show Corollaries \ref{c:gen_fun_sp} and
\ref{c:gen_fun_sl},
it suffices to consider a square matrix $X$ 
of the form 
$
X=\left[ \begin{smallmatrix}
          u 1_p & b \\ c & v 1_q
         \end{smallmatrix}
  \right]
$
with $u,v \in \C$ 
and submatrices $b$, $c$ of size 
$p \times q$, $q \times p$ respectively. 
Let us parametrize $X$ as follows:
\begin{equation*}
\label{e:tmp}
  X=\left[ \begin{array}{cccc|ccc}
     u    &       &       &       & b_{1,1} & \cdots & b_{1,q} \\
          & u     &       &       & b_{2,1} & \cdots & b_{2,q} \\
          &       &\ddots &       &  \vdots &        &  \vdots \\
          &       &       & u     & b_{p,1} & \cdots & b_{p,q} \\
    \hline
  c_{1,1} &c_{2,1}&\cdots &c_{p,1}& v       &        &         \\
  \vdots  &\vdots &       &\vdots &         &\ddots  &         \\
  c_{1,q} &c_{2,q}&\cdots &c_{p,q}&         &        & v
          \end{array}
   \right].
\end{equation*}

\begin{prop}
\label{p:msf_det}
Let $X$ be a square matrix given as above.
Then we can expand $\det X$
in the following way:
\begin{equation}
\label{e:msf_det}
  \det {X} = \sum_{k=0}^{q} 
          \sum_{ \substack{
                 I \subset [p], J \subset [q] \\%
                 |I|=|J|=k}  }
            u^{p-k} v^{q-k} \det{(b^{I}_{J})} \det{(c^{I}_{J})}.
\end{equation}
 
\end{prop}


\bibliographystyle{amsalpha}
\bibliography{rep}

\nocite{HOOS93}
\nocite{Shimura84}

\end{document}

%% file: env.tex
\usepackage[T1]{fontenc}
\usepackage{textcomp}
\usepackage{amsmath,amsthm,amssymb,amscd}
\usepackage{type1cm}
\usepackage{txfonts}

\numberwithin{equation}{section}	

\theoremstyle{plain}
\newtheorem{thm}{Theorem}[section]      
\newtheorem*{thm*}{Theorem}

\newtheorem{prop}[thm]{Proposition}

\newtheorem{cor}[thm]{Corollary}

\theoremstyle{definition}

\theoremstyle{remark}
\newtheorem{remark}[thm]{Remark}

\theoremstyle{remark}


%% file: ncpfaff_def.tex
\newcommand{\hph}[1]{\hphantom{#1}}

\def\geqsl{\geqslant}
\def\leqsl{\leqslant}
\newcommand{\floor}[1]{\lfloor{#1}\rfloor}
\newcommand{\tild}[1]{\widetilde{#1}}
\newcommand{\trace}[1]{\ensuremath{\operatorname{tr}{\left( #1 \right)}}}
\newcommand{\sgn}[1]{\ensuremath{\operatorname{sgn}{\left( #1 \right)}}}
\newcommand{\tp}[1]{\ensuremath{\operatorname{}^{t}\hspace{-3truept}{#1}}}

\newcommand{\Mat}[1]{\ensuremath{\operatorname{Mat}_{#1}}}
\newcommand{\Alt}{\ensuremath{\operatorname{Alt}}}   
   
\newcommand{\bs}[1]{\ensuremath{\boldsymbol{#1}}}
\newcommand{\Pf}[1]{\ensuremath{\operatorname{Pf}\left({#1}\right)}}

\newcommand{\Lie}[1]{\ensuremath{\operatorname{Lie}{#1}}}

\newcommand{\Rrank}{\ensuremath{\R\textrm{-}\operatorname{rank}}}

\newcommand{\Ad}{\ensuremath{\operatorname{Ad}}}  

\newcommand{\umat}[1]{E^{(#1)}}

\def\c{\gamma}

\def\H{\mathscr H}	%

\renewcommand{\Omega}{\ensuremath{\varOmega}}
\renewcommand{\Xi}{\ensuremath{\varXi}}
\renewcommand{\Theta}{\ensuremath{\varTheta}}

\def\R{\mathbb R}
\def\C{\mathbb C}

\newcommand{\gl}{\ensuremath{\mathfrak {gl}}}
\renewcommand{\sl}{\ensuremath{\mathfrak {sl}}}

\newcommand{\so}{\ensuremath{\mathfrak{so}}}
\renewcommand{\sp}{\ensuremath{\mathfrak{sp}}}

\newcommand{\g}{\ensuremath{\mathfrak{g}}}   
\newcommand{\p}{\ensuremath{\mathfrak{p}}}   
   
\renewcommand{\u}{\ensuremath{\mathfrak{u}}}

\newcommand{\GL}[1]{\ensuremath{\mathrm{GL}_{#1}}}
\newcommand{\SL}[1]{\ensuremath{\mathrm{SL}_{#1}}}
   
\newcommand{\SO}{\ensuremath{\mathrm{SO}}}   
\newcommand{\U}{\ensuremath{\mathrm{U}}}   
\newcommand{\SU}{\ensuremath{\mathrm{SU}}}   
\newcommand{\Sp}{\ensuremath{\mathrm{Sp}}}   
\newcommand{\SOstar}{\ensuremath{\mathrm{SO}^*}}   

\newcommand{\Gc}{\ensuremath{G_{\C}}}
\newcommand{\Kc}{\ensuremath{K_{\C}}}

\def\Z{\mathbb Z}

\renewcommand{\d}{\ensuremath{\operatorname{d}\!}}
\newcommand{\pd}{\ensuremath{\partial}}

\def\<{\langle}
\def\>{\rangle}

\def\c2vec#1#2{ %
   \left[ \begin{smallmatrix} %
           #1 \\ #2  \end{smallmatrix} %
   \right]}
